\newcommand\CC{{\mathbb C}}
\newcommand\cA{{\cal A}}

\newcommand\cB{{\cal B}}

\newcommand\cE{{\cal E}}
\newcommand\cF{{\cal F}}
\newcommand\cG{{\cal G}}
\newcommand\cH{{\cal H}}

\newcommand\cK{{\cal K}}
\newcommand\cL{{\cal L}}

\newcommand\cO{{\cal O}}
\newcommand\cP{{\cal P}}

\newcommand\cS{{\cal S}}

\newcommand\cU{{\cal U}}
\newcommand\cV{{\cal V}}
\newcommand\cW{{\cal W}}

\newcommand\hra{\hookrightarrow}
\newcommand\lra{\longrightarrow}

\newcommand\ov{\overline}
\newcommand\PP{{\mathbb P}}

\newcommand\wt{\widetilde}
\newcommand\ZZ{{\mathbb Z}}

\newcommand{\mapor}[1]{{\stackrel{#1}{\longrightarrow}}}
\newcommand{\mapver}[1]{\Big\downarrow
\vcenter{\rlap{$\scriptstyle#1$}}}
\documentclass[10pt,a4paper]{article}
\usepackage{amsthm}
\usepackage{amsmath}
\usepackage{amssymb}
\newtheorem{thm}{Theorem}[section]
\newtheorem{clm}[thm]{Claim}

\newtheorem{crl}[thm]{Corollary}
\newtheorem{dfn}[thm]{Definition}
 \newtheorem{lmm}[thm]{Lemma}
\newtheorem{prp}[thm]{Proposition}

\newtheorem{rmk}[thm]{Remark}

\begin{document}
 \title{Irreducible  symplectic $4$-folds and
 Eisenbud-Popescu-Walter sextics}
 \author{Kieran G. O'Grady\thanks{Supported by
 Cofinanziamento MIUR 2004-2005}\\
Universit\`a di Roma \lq\lq La Sapienza\rq\rq}
\date{July 19 2005
\vskip 1cm
\begin{abstract}
 Eisenbud Popescu and Walter have constructed
 certain special sextic
 hypersurfaces in $\PP^5$ as Lagrangian degeneracy loci. We
 prove that the natural double cover of a generic
 EPW-sextic is a deformation of the Hilbert
 square of a $K3$ surface $(K3)^{[2]}$ and that
 the family of such varieties is locally complete
 for deformations that keep the hyperplane class
 of type $(1,1)$ - thus we get an example similar to
 that (discovered by Beauville and Donagi)
 of the Fano variety of lines on a cubic
$4$-fold. Conversely suppose that $X$ is a
numerical $(K3)^{[2]}$, that $H$ is an ample
divisor on $X$ of square $2$ for Beauville's
quadratic form and that the map
$X\dashrightarrow|H|^{\vee}$ is the composition
of the quotient $X\to Y$ for an anti-symplectic
involution on $X$ followed by an immersion
$Y\hra|H|^{\vee}$; then $Y$ is an EPW-sextic and
$X\to Y$ is the natural double cover. If a
conjecture on the behaviour of certain linear
systems holds this result together with previous
results of ours implies that every numerical
$(K3)^{[2]}$ is a deformation of $(K3)^{[2]}$.
 \end{abstract}}
 \maketitle
 \section{Introduction}\label{prologo}
 \setcounter{equation}{0}
A compact K\"ahler manifold is
 irreducible symplectic if it is
simply connected and it carries a holomorphic
symplectic form spanning $H^{2,0}$
(see~\cite{beau,huy}). An irreducible
symplectic surface is nothing else but a $K3$
surface. Higher-dimensional  irreducible
symplectic manifolds behave like $K3$ surfaces
in many respects~\cite{huy,huy-err} however
their classification up to deformation of
complex structure is out of reach at the
moment. Let $S$ be a $K3$; the Hilbert square
$S^{[2]}$ i.e.~the blow-up of the diagonal in
the symmetric square $S^{(2)}$ is the symplest
example of an irreducible symplectic $4$-fold.
An irreducible symplectic
 $4$-fold $M$ is a {\it numerical $(K3)^{[2]}$}
 if there exists an isomorphism
 of abelian groups
\begin{equation}\label{}
  \psi\colon
 H^2(M;\ZZ)\overset{\sim}{\lra}H^2(S^{[2]};\ZZ)
\end{equation}
 such that
\begin{equation}\label{}
\int_{M}\alpha^{4}=\int_{S^{[2]}}\psi(\alpha)^{4}
\end{equation}
for all $\alpha\in H^2(M;\ZZ)$\footnote{For the
experts: the Beauville quadratic form and the
Fujiki constant of $M$ are the same as those of
$S^{[2]}$.}. In~\cite{ognum} we studied
numerical $(K3)^{[2]}$'s with the goal of
classifying them up to deformation of complex
structure. We proved that any numerical
$(K3)^{[2]}$ is deformation equivalent to an
$X$ carrying an ample divisor $H$ such that
\begin{equation}\label{gradodod}
  \int_X c_1(H)^4=12,\quad \dim |H|=5
\end{equation}
(the first equation is equivalent to $c_1(H)$
being of square $2$ for the Beauville form),
and such that the rational map
\begin{equation}\label{mappaeffe}
  X\dashrightarrow |H|^{\vee}
\end{equation}
satisfies one of the following two conditions:
\begin{itemize}
  \item [(a)]
 There exist an
 anti-symplectic involution $\phi\colon X\to X$
 with quotient $Y$  and an embedding $Y\hra
 |H|^{\vee}$ such that~(\ref{mappaeffe}) is the
 composition of the quotient map $f\colon X\to
 Y$ and the embedding $Y\hra
 |H|^{\vee}$.
  \item [(b)]
Map~(\ref{mappaeffe}) is birational onto a
hypersurface of degree between $6$ and $12$.
\end{itemize}
 In this paper we describe
explicitely all the $X$ occuring in Item~(a)
above. Notice that $Y$ is singular because smooth
hypersurfaces in $\PP^5$ are simply connected.
Moreover the singular locus is a surface because
$\phi$ is anti-symplectic. Thus $Y$ is far from
being a generic sextic hypersurface; we will show
that it belongs to a family of sextics
constructed by Eisenbud, Popescu and Walter, see
Example~(9.3) of~\cite{epw}. We will prove that
conversely a generic EPW-sextic has a natural
double cover which is a deformation of
$(K3)^{[2]}$.  Since EPW-sextics form an
irreducible family we get that the $X$'s
satisfying~(a) above are deformation equivalent.
Actually if $(X_i,f_i^{*}\cO_{Y_i}(1))$ are two
polarized couples where $f_i\colon X_i\to Y_i$
satisfy~(a) above for $i=1,2$ then we may deform
$(X_1,f_1^{*}\cO_{Y_1}(1))$ to
$(X_2,f_2^{*}\cO_{Y_2}(1))$ through polarized
deformations. In particular all the explicit
examples of $f\colon X\to Y$ satisfying~(a) above
that were constructed in~\cite{oginv} are
equivalent through polarized deformations - this
answers positively a question raised in
Section~(6) of~\cite{oginv}. We recall that no
examples are known of $X$ satisfying Item~(b)
above; in~\cite{ognum} we conjectured that they
do not exist - one result in favour of the
conjecture is that
 if $X$ satisfies~(a) above then all small deformations
 of $X$ keeping $c_1(f^{*}\cO_Y(1))$
of type $(1,1)$ also satisfy~(a) above see the
proposition at the end of Section~(4)
of~\cite{ognum}. If our conjecture is true then
the results of this paper together with the
quoted results of~\cite{ognum} give that
numerical $(K3)^{[2]}$'s are deformation
equivalent to the Hilbert square of a $K3$.
Before stating precisely our main results we
recall the construction of EPW-sextics. Let $V$
be a $6$-dimensional vector space and $\PP(V)$ be
the projective space of $1$-dimensional sub
vector spaces $\ell\subset V$. Choose an
isomorphism $vol\colon\wedge^6
V\overset{\sim}{\lra}\CC$ and let $\sigma$ be the
symplectic form on $\wedge^3 V$ defined by wedge
product,
i.e.~$\sigma(\alpha,\beta):=vol(\alpha\wedge\beta)$;
thus $\wedge^3 V\otimes\cO_{\PP(V)}$ has the
structure of a symplectic vector-bundle of rank
$20$. Let $F$ be  the sub-vector-bundle of
$\wedge^3 V\otimes\cO_{\PP(V)}$ with fiber over
$\ell\in\PP(V)$ equal to
\begin{equation}\label{fibradieffe}
  F_{\ell}:=
  Im\left(\ell\otimes\wedge^2(V/\ell)
  \hra\wedge^3 V\right).
\end{equation}
Thus we have an exact sequence
\begin{equation}\label{effeinwedge}
  0\to F\overset{\nu}{\lra}\wedge^3
  V\otimes\cO_{\PP(V)}
\end{equation}
We have $rk(F)=10$ and $\sigma|_{F_{\ell}}=0$;
thus $F$ is a Lagrangian sub-bundle of
$\wedge^3 V\otimes\cO_{\PP(V)}$. Let
$\mathbb{LG}(\wedge^3 V)$ be the symplectic
Grassmannian parametrizing $\sigma$-Lagrangian
subspaces of $\wedge^3 V$. For $A\in
\mathbb{LG}(\wedge^3 V)$ we let
\begin{equation}\label{mappaquot}
  \lambda_A\colon F\lra (\wedge^3 V/A)\otimes\cO_{\PP(V)}
\end{equation}
be Inclusion~(\ref{effeinwedge}) followed by
the projection $(\wedge^3
V)\otimes\cO_{\PP(V)}\to (\wedge^3
V/A)\otimes\cO_{\PP(V)}$. Since the
vector-bundles appearing in~(\ref{mappaquot})
have equal rank we have $\det(\lambda_A)\in
H^0((\det F)^{-1})$. We let $Y_A\subset\PP(V)$
be the zero-scheme of $\det(\lambda_A)$. Let
$\omega:=c_1(\cO_{\PP(V)}(1))$; a
straightforward computation gives that
\begin{equation}\label{cherndieffe}
  c(F)=1-6\omega+18\omega^2-34\omega^3+\ldots
\end{equation}
In particular $\det F\cong\cO_{\PP(V)}(-6)$.
Thus $Y$ is always non-empty and if
$Y\not=\PP(V)$ then $Y$ is a sextic
hypersurface. An {\it EPW-sextic} is a
hypersurface in $\PP(V)$ which is equal to
$Y_A$ for some $A\in \mathbb{LG}(\wedge^3 V)$.
In Section~(\ref{eccoepw}) we describe
explicitely the non-empty Zariski-open
$\mathbb{LG}(\wedge^3
V)^0\subset\mathbb{LG}(\wedge^3 V)$
parametrizing $A$ such that the following hold:
$Y_A$ is a sextic hypersurface smooth at all
points where the map $\lambda_A$
of~(\ref{mappaquot}) has corank $1$, the
analytic germ $(Y_A,\ell)$ at a point $\ell$
where $\lambda_A$ has corank $2$
 is isomorphic to the product of a
smooth $2$-dimensional germ times an
$A_1$-singularity and furthermore $\lambda_A$
 has corank at most $2$ at all points of
 $\PP(V)$. Let
 $A\in\mathbb{LG}(\wedge^3 V)^0$; then
  $Y_A$ supports a
 quadratic sheaf as defined by Casnati and
 Catanese~\cite{cascat} and hence
 there is a natural double cover $X_A\to Y_A$
 with $X_A$ smooth - see Section~(\ref{epwdoppio}). In
 Section~(\ref{dimprinc}) we will prove the
 following result.
 \begin{thm}\label{mainthm}
 Keep notation as above. Then the following
 hold.
\begin{itemize}
  \item [(1)]
 Suppose that $X,H$ are a
 numerical $(K3)^{[2]}$ and an ample divisor on
 $X$ such that both~(\ref{gradodod}) and Item~(a)
 above hold. Then there exists
 $A\in\mathbb{LG}(\wedge^3
 V)^0$ such that $f\colon X\to Y$ is identified
 with the natural double cover $X_A\to
 Y_A$.
  \item [(2)]
 For $A\in\mathbb{LG}(\wedge^3 V)^0$
 let $X_A\to Y_A$ be the natural double cover
defined in Section~(\ref{epwdoppio}) and let
$H_A$ be the pull-back to $X_A$ of
$\cO_{Y_A}(1)$. Then $X_A$ is an irreducible
symplectic variety deformation equivalent to
$(K3)^{[2]}$ and both~(\ref{gradodod}) and
Item~(a) above hold with $X=X_A$ and $H=H_A$.
\end{itemize}
 \end{thm}
We recall that Beauville-Donagi~\cite{beaudon}
proved the following result: if $Z\subset\PP^5$
is a smooth cubic hypersurface the Fano variety
$F(Z)$ parametrizing lines on $Z$ is a
deformation of $(K3)^{[2]}$. They also proved
that the family of $F(Z)$'s polarized
 by the Pl\"ucker
line-bundle is locally complete. Similarly the family of $X_A$'s that one gets
by letting $A$ vary in
 $\mathbb{LG}(\wedge^3 V)^0$ is also a locally
 complete family of polarized varieties by the
 proposition at the end
  of Section~(4) of~\cite{ognum} that we have
  already quoted. This is confirmed by the
  following computation. The tangent
   space to $\mathbb{LG}(\wedge^3 V)$ at a
   point $A$ is isomorphic to $Sym_2 A^{\vee}$
 and hence $\dim\mathbb{LG}(\wedge^3 V)=55$. Since
  $\dim\mathbb{PGL}(V)=35$ we get that
  $\dim(\mathbb{LG}(\wedge^3 V)^0//
 \mathbb{PGL}(V))=20$ which is the number of
 moduli of a polarized deformation of
 $(K3)^{[2]}$. We remark that the Beauville-Donagi family
 and the family of $X_A$'s
 are the only explicit examples of a locally complete family
 of higher dimensional polarized
 irreducible symplectic varieties. Notice that our
 conjecture amounts to the statement that the familiy
  of $X_A$'s is
 globally complete once we take into account
 the limiting $X_A$'s one gets for
 $A\in\left(\mathbb{LG}(\wedge^3 V)\setminus
 \mathbb{LG}(\wedge^3 V)^0\right)$. An interesting
 feature of EPW-sextics is that
they are preserved by the duality map i.e.~the
dual of an EPW-sextic $Y_A$ is an EPW-sextic,
see Section~(\ref{dualediepw}). Thus duality
defines a regular involution on an open dense
subset of the moduli space of numerical
$(K3)^{[2]}$'s polarized by a divisor $H$
satisfying~(\ref{gradodod})\footnote{The second
equation of~(\ref{gradodod}) follows from the
first equation - see~\cite{ognum}.} and
Item~(a) above, see Section~(\ref{involsuk}).
It would be interesting to know explicitly
which $A\in \mathbb{LG}(\wedge^3
 V)$ correspond to  special $4$-folds
  e.g.~Hilbert squares of a $K3$; we discuss
 this problem in Section~(\ref{esespl}).
 \vskip 3mm

 \noindent
 {\bf Acknowledgement:} It is a pleasure to
 thank Adrian Langer for the interest he took in
 this work. In
 particular Adrian proved
 Proposition~(\ref{difftre}) and indicated how to
 prove Proposition~(\ref{treacca}).
 \section{EPW-sextics}
 \label{eccoepw}
 \setcounter{equation}{0}
  We will explicitely describe
   those EPW-sextics whose only
  singularities are the expected ones - the main
  result is Proposition~(\ref{yewgen}).
 We start by recalling (see~\cite{epw,fulpra})
 how one defines natural
   subschemes $D_i(A,F)\subset\PP(V)$    such that
\begin{equation}\label{supportodid}
 supp D_i(A,F)=\{\ell\in\PP(V)|\ \
\dim(F_{\ell}\cap A)\ge i\}.
\end{equation}
\begin{dfn}\label{banalsimpl}
Let $U\subset\PP(V)$ be an open
 subset. A {\it symplectic trivialization of $\wedge^3
 V\otimes\cO_U$} consists of a couple $(\cL,\cH)$ of trivial transversal Lagrangian
 sub-vector-bundles $\cL,\cH\subset\wedge^3
 V\otimes\cO_U$.
\end{dfn}
Let $(\cL,\cH)$ be a symplectic trivialization of
$\wedge^3
 V\otimes\cO_U$. The symplectic form $\sigma$ defines
 an isomorphism
\begin{equation}\label{}
  \begin{matrix}
 \cH & \overset{\sim}{\lra} & \cL^{\vee} \\
 \alpha & \mapsto & \sigma(\alpha,\cdot)
  \end{matrix}
\end{equation}
  and hence we get a direct sum
 decomposition
\begin{equation}\label{elpiuelduale}
\cL\oplus\cL^{\vee}= \wedge^3 V\otimes\cO_U.
\end{equation}
Conversely a direct sum
decomposition~(\ref{elpiuelduale}) with
$\cL,\cL^{\vee}$ trivial Lagrangian
sub-vector-bundles such that $\sigma$ induces
the tautological isomorphism
$\cL^{\vee}\overset{\sim}{\lra}\cL^{\vee}$
gives
 a symplectic trivialization of $\wedge^3
V\otimes\cO_U$; this is how we usually present a
symplectic trivialization of $\wedge^3
V\otimes\cO_U$.
 \begin{clm}\label{buonabanana}
 Let $A\in\mathbb{LG}(\wedge^3 V)$
 and $\ell_0\in\PP(V)$. There exists an open affine
 $U\subset\PP(V)$ containing $\ell_0$ and a
 symplectic trivialization~(\ref{elpiuelduale}) of
 $\wedge^3 V\otimes\cO_U$ such that for every $\ell\in U$
 both $A$ and $F_{\ell}$ are  transversal to
 $\cL^{\vee}_{\ell}$.
\end{clm}
\begin{proof}
The set of Lagrangian subspaces of $(\wedge^3 V)$
which are transversal to a given Lagrangian
subspace is an open dense subset of
$\mathbb{LG}(\wedge^3 V)$. Thus there exists
$C\in\mathbb{LG}(\wedge^3 V)$ which is
transversal both to $A$ and to $F_{\ell_0}$.
Since the condition of being transversal is open
there exists an open affine $U\subset\PP(V)$
containing $\ell_0$ such that $F_{\ell}$ is
transversal to $C$ for all $\ell\in U$.
 Let $B\in\mathbb{LG}(\wedge^3 V)$ be
 transversal to $C$; thus
\begin{equation}\label{cipiudi}
  \wedge^3 V=B\oplus C
\end{equation}
and hence $(B\otimes\cO_U,C\otimes\cO_U)$ is a
symplectic trivialization of $\wedge^3
V\otimes\cO_U$. Letting $\cL:=B\otimes\cO_U$ we
have an isomorphism $\cL^{\vee}\cong
C\otimes\cO_U$ induced by $\sigma$ and we write
the chosen symplectic trivialization
as~(\ref{elpiuelduale}); by construction both
$A$ and $F_{\ell}$ are transversal to
$\cL^{\vee}_{\ell}$ for every $\ell\in U$.
\end{proof}
Choose $A\in\mathbb{LG}(\wedge^3 V)$.
 Let $\ell_0\in\PP(V)$. Assume that we have
$U$ and a symplectic trivialization of $\wedge^3
V\otimes\cO_U$ as in the claim above. We will
define a closed degeneracy subscheme
$D_i(A,F,U,\cL,\cL^{\vee})\subset U$ - after
doing this we will define the subcheme
$D_i(A,F)\subset\PP(V)$ by gluing together the
local degeneracy loci. Via~(\ref{elpiuelduale})
we may identify both $A\otimes\cO_U$ and $F|_U$
with the graphs of maps
\begin{equation}\label{mappegrafi}
  q_A,q_F\colon \cL\to\cL^{\vee}
\end{equation}
because for every $\ell\in U$ both $A$ and
$F_{\ell}$ are  transversal to
$\cL^{\vee}_{\ell}$. Since both $A\otimes\cO_U$
and $F|_U$ are Lagrangian sub-vector-bundles  the
maps $q_A$ and $q_F$ are symmetric. Choosing a
trivialization of $\cL$ we view $q_A,q_F$ as
symmetric $10\times 10$ matrices with entries in
$\CC[U]$.
\begin{dfn}\label{eccoluogodeg}
Keep notation as above. We let $D_i(A,F,U,\cL,\cL^{\vee})\subset U$ be the
closed subscheme defined by the vanishing of determinants of
$(11-i)\times(11-i)$-minors of $(q_A-q_F)$.
\end{dfn}
Notice that the definition above makes sense because if we change the
trivialization of $\cL$ the relevant determinants are multiplied by units of
$\CC[U]$
\begin{lmm}\label{ugualisuint}
Let $U_1,U_2\subset\PP(V)$ be open affine and
$\cL_j\oplus\cL^{\vee}_j=\wedge^3
V\otimes\cO_{U_j}$ be symplectic trivializations.
Then
\begin{equation}\label{}
  D_i(A,F,U_1,\cL_1,\cL_1^{\vee})\cap U_1\cap
U_2=D_i(A,F,U_2,\cL_2,\cL_2^{\vee})\cap U_1\cap U_2.
\end{equation}
\end{lmm}
\begin{proof}
It suffices to prove the lemma for $U_1=U_2=U$.
The constructions above can be carried out more
generally for a trivialization
$\cV\oplus\cW\cong\wedge^3
 V\otimes\cO_U$ where $\cV,\cW\subset\wedge^3
 V\otimes\cO_U$ are trivial rank-$10$ sub-vector-bundles
(not necessarily Lagrangian) such that  for every $\ell\in U$ both $A$ and
$F_{\ell}$ are transversal to $\cW_{\ell}$. We identify $A\otimes\cO_U$ and
$F|_U$ with the graphs of maps $q_A\colon \cV\to\cW$ and $q_F\colon \cV\to\cW$
respectively - notice that in general it does not make sense to ask whether
$q_A$, $q_F$ are symmetric! Trivializing $\cV$ and $\cW$ we view $q_A,q_F$ as
$10\times 10$ matrices with entries in $\CC[U]$. We let
$D_i(A,F,U,\cV,\cW)\subset U$ be the subscheme defined by the vanishing of
determinants of $(11-i)\times(11-i)$-minors of $(q_A-q_F)$. One checks easily
that if we change $\cV$ (leaving $\cW$ fixed) or if we change $\cW$ (leaving
$\cV$ fixed) the scheme $D_i(A,F,U,\cV,\cW)$ remains the same. The lemma follows
immediately.
\end{proof}
Now we define $D_i(A,F)$. Consider the collection
of symplectic trivializations
$\cL_j\oplus\cL_j^{\vee}=\wedge^3
V\otimes\cO_{U_j}$ with $U_j\subset\PP(V)$ open
affine and the corresponding closed subschemes
$D_i(A,F,U_j,\cL_j,\cL_j^{\vee})$. By
Claim~(\ref{buonabanana}) the subsets $U_j$ cover
$\PP(V)$ and by Lemma~(\ref{ugualisuint}) the
$D_i(A,F,U_j,\cL_j,\cL_j^{\vee})$ for different
$j$'s match on overlaps; thus they glue together
and they define a closed subscheme
$D_i(A,F)\subset\PP(V)$.
Clearly~(\ref{supportodid}) holds and furthermore
$D_{i+1}(A,F)$ is a subscheme of $D_i(A,F)$. We
claim that
\begin{equation}\label{}
 Y_A:=D_1(A,F).
\end{equation}
It is clear from~(\ref{supportodid}) that
$supp(Y_A)=supp(D_1(A,F))$ and hence we need
only check that the scheme structures coincide
in a neighborhood of any point $\ell_0\in
D_1(A,F)$. There exists $B\in{\mathbb
LG}(\wedge^3 V)$ which is transversal both to
$F_{\ell_0}$ and $A$. There is an open
neighborhood $U$ of $\ell_0$ such that $B$ is
transversal to $F_{\ell}$ for $\ell\in U$. The
symplectic form $\sigma$ defines an isomorphism
$B\cong A^{\vee}$. Consider the symplectic
trivialization of $\wedge^3V\otimes\cO_U$ given
by $\cL:=A\otimes\cO_{U}$ and
$\cL^{\vee}:=B\otimes\cO_{U}$; since $q_A=0$ we
have $D_1(A,F,U,\cL,\cL^{\vee})=Y_A\cap U$ and
we are done. To simplify notation we let
\begin{equation}\label{}
  W_A:=D_2(A,F).
\end{equation}
As shown in Section~(\ref{prologo}) - see~(\ref{mappaquot}) - $Y_A$ is never
empty. We claim that $W_A$ is never empty as well. In fact
 Formula~(6.7)
of~\cite{fulpra} and Equation~(\ref{cherndieffe})
give that if $W_A$ has the expected
 dimension i.e.~$2$ (see
 Equations~(\ref{diconipull})-(\ref{codimsigma})) or
 is empty then the cohomology class of the cycle
 $[W_A]$ is
\begin{equation}\label{classew}
  cl([W_A])=2c_3(F)-c_1(F)c_2(F)=40\omega^3.
\end{equation}
Since the right-hand side of the above equation
is non-zero it follows that necessarily
$W_A\not=\emptyset$.
\begin{dfn}\label{eccograsszero}
Let $\mathbb{LG}(\wedge^3 V)^{\times}\subset
\mathbb{LG}(\wedge^3 V)$ be the set of $A$ such
that for all $\ell\in\PP(V)$ we have
\begin{equation}\label{alpiudue}
  \dim A\cap F_{\ell}\le 2,
\end{equation}
i.e.~$D_3(A,F)=\emptyset$. Let
$\mathbb{LG}(\wedge^3 V)^0\subset
\mathbb{LG}(\wedge^3 V)^{\times}$ be the set of
$A$ which do not contain a non-zero completely
decomposable element $v_0\wedge v_1\wedge v_2$.
\end{dfn}
A straightforward dimension count gives the following
result.
\begin{clm}
Both $\mathbb{LG}(\wedge^3 V)^{\times}$ and
$\mathbb{LG}(\wedge^3 V)^0$ are open dense
subsets of $\mathbb{LG}(\wedge^3 V)$.
\end{clm}
We will show that $\mathbb{LG}(\wedge^3 V)^0$ is
the open subset of $\mathbb{LG}(\wedge^3 V)$
parametrizing $A$ such that $Y_A$ and $W_A$ are
as nice as possible. First we describe
$\Theta_{\ell_0}D_i(A,F)$ at a point $\ell_0\in
D_i(A,F)$. Proceeding as in the proof of
Claim~(\ref{buonabanana}) we consider a
symplectic trivialization
\begin{equation}\label{}
  \wedge^3 V\otimes\cO_U=A\otimes\cO_U\oplus
  A^{\vee}\otimes\cO_U
\end{equation}
where $U\subset\PP(V)$ is a suitable open affine
subset containing $\ell_0$ and
$A^{\vee}\in\mathbb{LG}(\wedge^3 V)$ is
transversal  to $A$ and to $F_{\ell}$ for every
$\ell\in U$.  We view $F\otimes\cO_U$ as the
graph of a symmetric map $q_F\colon
A\otimes\cO_U\to A^{\vee}\otimes\cO_U$. Let
\begin{equation}\label{mappapsi}
   \begin{matrix}
 U & \overset{\psi}{\lra} & Sym_2 A^{\vee} \\
 \ell & \mapsto & \psi(\ell):=q_F(\ell)
   \end{matrix}
\end{equation}
and let $\Sigma_i\subset Sym_2 A^{\vee}$ be the
closed subscheme parametrizing quadratic forms of
corank at least $i$. Since the map $q_A\colon
A\otimes\cO_U\to A^{\vee}\otimes\cO_U$ whose
graph is $A\otimes\cO_U$ is zero we have
\begin{equation}\label{diconipull}
  D_i(A,F)\cap U=\psi^{*}\Sigma_i .
\end{equation}
 We recall that $\Sigma_i$ is an irreducible
 local complete intersection with
\begin{equation}\label{codimsigma}
  cod(\Sigma_i,Sym_2 A^{\vee})=i(i+1)/2.
\end{equation}
Let $\ov{q}\in(\Sigma_i\setminus\Sigma_{i+1})$. Then $\Sigma_i$ is smooth at
$\ov{q}$ and the tangent space $\Theta_{\ov{q}}\Sigma_i$  is described as
follows. Identify $Sym_2 A^{\vee}$  with its tangent space at $\ov{q}$; then
\begin{equation}\label{conotangente}
 \Theta_{\ov{q}}\Sigma_i=\{q\in Sym_2 A^{\vee}|\
 q|_{\ker(\ov{q})}=0\}.
\end{equation}
Thus we also get a natural identification
\begin{equation}\label{spazionorm}
 (N_{\Sigma_i/Sym_2 A^{\vee}})_{\ov{q}}=
 Sym_2 \ker(\ov{q})^{\vee}.
\end{equation}
\begin{lmm}\label{faine}
Keep notation as above. Suppose that
 $A\in\mathbb{LG}(\wedge^3 V)^{\times}$.
Let $\ell_0\in D_i(A,F)$ and let
\begin{equation}\label{}
  A\cap F_{\ell_0}=\ell_0\otimes W
\end{equation}
where $W\subset \wedge^2(V/\ell_0)$. The
composition
\begin{equation}\label{ranieri}
  \Theta_{\ell_0}\PP(V)\overset{d\psi(\ell_0)}{\lra}
  Sym_2 A^{\vee}\lra Sym_2 \ker\psi(\ell_0)^{\vee}
\end{equation}
is surjective if and only if $W$ contains no non-zero
decomposable element.
\end{lmm}
\begin{proof}
 Let $\ell_0=\CC v_0$. Choose a codimension $1$
 subspace $V_0\subset V$
 transversal to $\ell_0$. Thus we have
\begin{equation}\label{}
  V/\ell_0\cong\oplus V_0,\quad
 W\subset\wedge^2 V_0.
\end{equation}
The map
\begin{equation}\label{}
  \begin{matrix}
 V_0 & \lra & \PP(V) \\
 u & \mapsto & [v_0+u]
  \end{matrix}
\end{equation}
gives an isomorphism between $V_0$
 and an open affine subspace of $\PP(V)$ containing
 $\ell_0$ - with $0\in V_0$ corresponding to $\ell_0$.
 Shrinking $U$ if necessary we may assume that
  $U\subset V_0$ is an open subset containing
  $0$. For $u\in U$ the map $\psi(u)\colon A\to
  A^{\vee}$ is characterized by the  equation
\begin{equation}\label{}
 \alpha+\psi(u)(\alpha)=(v_0+u)\wedge\gamma(u,\alpha),
 \quad \gamma(u,\alpha)\in\wedge^2 V_0
\end{equation}
 where $\alpha\in A$. Thus when we view $\psi(u)$
 as a symmetric bilinear form we have the formula
 \begin{equation}\label{}
 \psi(u)(\alpha,\beta)=
 vol\left((v_0+u)\wedge\gamma(u,\alpha)\wedge\beta\right)
\end{equation}
 for $\alpha,\beta\in A$. Now assume that
 $\alpha,\beta\in\ker\psi(\ell_0)$ and hence
\begin{equation}\label{}
  \alpha=v_0\wedge\alpha_0,\quad
  \beta=v_0\wedge\beta_0,\quad
  \alpha_0,\beta_0\in\wedge^2 V_0.
\end{equation}
 Let
\begin{equation}\label{}
  \tau\in\Theta_0 U\cong V_0
\end{equation}
and let $u(t)$ be a
 \lq\lq parametrized curve\rq\rq in $U$ with
 $u(0)=0$ and $\dot{u}(0)=\tau$. Then
 \begin{equation}\label{}
 d\psi(\tau)(v_0\wedge\alpha_0,
 v_0\wedge\beta_0)=
 \frac{d}{dt}_{\vert t=0}
 vol\left((v_0+u(t))\wedge\gamma(u(t),v_0\wedge\alpha_0)
 \wedge v_0\wedge\beta_0\right).
\end{equation}
 Differentiating and observing that
 $\gamma(0,v_0\wedge\alpha_0)=\alpha_0$ we get
 that
\begin{equation}\label{eccoderivata}
 d\psi(\tau)(v_0\wedge\alpha_0,
 v_0\wedge\beta_0)=-vol\left(v_0\wedge
 \tau\wedge\alpha_0\wedge\beta_0\right).
\end{equation}
Let's prove the proposition. If $i=0$ there is
nothing to prove. If $i=1$ let
$\ker\psi(\ell_0)=\CC v_0\wedge\alpha_0$. We
apply the above formula with $\beta_0=\alpha_0$;
since $\tau\in V_0$ is an arbitrary element
 we get that
Composition~(\ref{ranieri}) is surjective if and
only if $\alpha_0\wedge\alpha_0\not=0$ i.e.~if
and only if $\alpha_0$ is not decomposbale. This
proves the proposition when $i=1$. By definition
of $\mathbb{LG}(\wedge^3 V)^{\times}$ we are left
with the case $i=2$. First  notice that if
$\alpha_0\in W$ is decomposable then
$d\psi(\tau)(v_0\wedge\alpha_0,
 v_0\wedge\alpha_0)=0$ for all $\tau$ by
 Equation~(\ref{eccoderivata}); thus if
 $W$ contains a non-zero
 decomposable element then
 Composition~(\ref{ranieri}) is not surjective.
 Now assume that $W$ does not contain
  non-zero decomposable elements. Then we have a
  well-defined regular map
\begin{equation}\label{}
   \begin{matrix}
 \PP(W) & \overset{\rho}{\lra} & \PP(\wedge^4 V_0) \\
 [\alpha] & \mapsto & [\alpha\wedge\alpha]
   \end{matrix}
\end{equation}
We claim that $\rho$ is injective. In fact assume
that we have $[\alpha]\not=[\beta]$ and
$[\alpha\wedge\alpha]=[\beta\wedge\beta]$. Then
$span(\alpha)=span(\beta)=S$ where $S\subset V_0$
is a subspace of dimension $4$. Since
$\alpha,\beta$ span $W$ we get that
$span(\gamma)=S$ for all $\gamma\in W$. Thus
\begin{equation}\label{}
  \gamma\wedge\gamma\in
  Im(\wedge^2 S\to \wedge^2 V_0)
\end{equation}
for all $\gamma\in W$. Since
$\mathbb{G}r(2,S)\subset\PP(\wedge^2 S)$ is a
hypersurface we get that there exists a
decomposable non-zero $\gamma\in W$,
contradiction. Thus $\rho$ is injective; since
$\rho$ is defined by quadratic polynomials we get
that $Im(\rho)$ is a conic in $\PP(\wedge^4 V_0)$
and hence
\begin{equation}\label{}
 \wedge^4 V_0^{\vee}\lra H^0(\cO_{\PP(W)}(2))
\end{equation}
is surjective. Given Formula~(\ref{eccoderivata})
this implies surjectivity of
Composition~(\ref{ranieri}).
\end{proof}
\begin{prp}\label{yewgen}
Keep notation as above. Suppose that $A\in
\mathbb{LG}(\wedge^3 V)^{\times}$. The
following statements are equivalent.
\begin{itemize}
  \item [(1)]
 $A\in \mathbb{LG}(\wedge^3 V)^0$.
  \item [(2)]
$(D_i(A,F)\setminus D_{i+1}(A,F))$ is smooth of
codimension equal to the expected codimension
$i(i+1)/2$ for all $i$.
  \item [(3)]
 $(Y_A\setminus W_A)$ is smooth and for every
 $\ell_0\in W_A$ the following holds. There exist
 $U\subset\PP(V)$ open in the analytic topology
  containing $\ell_0$ and $x,y,z$ holomorphic
  functions on $U$ vanishing at $\ell_0$ with
  $dx(\ell_0),dy(\ell_0),dz(\ell_0)$ linearly
  independent such that
  $Y_A\cap U=V(xz-y^2)$.
\end{itemize}
\end{prp}
\begin{proof}
We prove equivalence of~(1) and~(2) - the proof
of equivalence of~(1) and~(3) is similar, we
leave it to the reader. Let $\ell_0\in
(D_i(A,F)\setminus D_{i+1}(A,F))$.
 Since we have Equation~(\ref{diconipull}) and
 since $(\Sigma_i\setminus\Sigma_{i+1})$ is smooth
 of codimension $i(i+1)/2$ (see
 Equation~(\ref{codimsigma})) $D_i(A,F)$ is
 smooth of codimension $i(i+1)/2$ at $\ell_0$ if
 and only if Composition~(\ref{ranieri}) is
 surjective. Let $\ell_0=\CC v_0$. By
 Lemma~(\ref{faine}) Composition~(\ref{ranieri}) is
 surjective if and only if $A$ does not contain a
 non-zero completely decomposable element
 divisible by $v_0$, i.e.~of the form $v_0\wedge
 v_1\wedge v_3$. The proposition follows immediately.
\end{proof}
\begin{rmk}\label{modelloloc}
 {\rm Let $A\in \mathbb{LG}(\wedge^3 V)^0$.
\begin{itemize}
  \item [(1)]
  By Proposition~(\ref{yewgen}) $Y_A$ has canonical
 singularities; since $Y_A$ is a sextic
 adjunction gives that $Y_A$
 has Kodaira dimension $0$.
  \item [(2)]
 The family of
 $Y_A$'s (for $A\in \mathbb{LG}(\wedge^3 V)^0$) is
 locally trivial.
\end{itemize}
}
 \end{rmk}
 \section{The dual of an EPW-sextic}
 \label{dualediepw}
 \setcounter{equation}{0}
Let $A\in\mathbb{LG}(\wedge^3 V)$ be  such that
$Y_A$ is a reduced hypersurface. For $\ell\in
Y_A^{sm}$ a smooth point of $Y_A$ let
$T_{\ell}Y_A\subset\PP(V^{\vee})$ be the
projective tangent space to $Y_A$ at $\ell$;
the dual $Y_A^{\vee}\subset\PP(V^{\vee})$ is
(as usual) the closure of $\bigcup_{\ell\in
Y^{sm}_A}T_{\ell}Y_A$. We will  show that if
$A$ is generic then $Y^{\vee}_A$  is isomorphic
to an EPW-sextic. The trivialization $vol\colon
\wedge^6 V\overset{\sim}{\lra}\CC$ defines a
trivialization $vol^{\vee}\colon \wedge^6
V^{\vee}\overset{\sim}{\lra}\CC$ and hence a
symplectic form $\sigma^{\vee}$ on $\wedge^3
V^{\vee}$; let $\mathbb{LG}(\wedge^3 V^{\vee})$
be the symplectic Grassmannian parametrizing
$\sigma^{\vee}$-Lagrangian subspaces of
$\wedge^3 V^{\vee}$. For
$A\in\mathbb{LG}(\wedge^3 V)$ we let
\begin{equation}\label{}
  A^{\bot}:=\{\phi\in \wedge^3 V^{\vee}|\
  \langle\phi,A\rangle=0\}
  \subset\wedge^3 V^{\vee}.
\end{equation}
As is easily checked $A^{\bot}\in
\mathbb{LG}(\wedge^3 V^{\vee})$. Thus we have
an isomorphism
\begin{equation}\label{annichilo}
  \begin{matrix}
 \delta_V\colon\mathbb{LG}(\wedge^3 V) &
 \overset{\sim}{\lra} &
 \mathbb{LG}(\wedge^3 V^{\vee})\\
 A & \mapsto & A^{\bot}.
  \end{matrix}
\end{equation}
\begin{prp}\label{dualediy}
 Keep notation as above and assume that
\begin{equation}\label{pigna}
  A\in\mathbb{LG}(\wedge^3
 V)^0\cap\delta_V^{-1}\mathbb{LG}(\wedge^3
 V^{\vee})^0.
\end{equation}
Then $Y^{\vee}_A=Y_{A^{\bot}}$.
\end{prp}
 \begin{proof}
 We claim that it suffices to prove that
\begin{equation}\label{contenuto}
  Y_A^{\vee}\subset Y_{A^{\bot}}.
\end{equation}
 In fact
 by Item~(1) of Remark~(\ref{modelloloc}) we know that
 $Y_A$ is not covered by positive-dimensional
 linear spaces and hence $Y_A^{\vee}$ is
 $4$-dimensional; since $A^{\bot}\in
 \mathbb{LG}(\wedge^3 V^{\vee})^0$ we know that
 $Y_{A^{\bot}}$ is irreducible (see
 Proposition~(\ref{yewgen}))  and
 hence~(\ref{contenuto}) implies that
 $Y^{\vee}_A=Y_{A^{\bot}}$.
 Let's prove~(\ref{contenuto}). Let $\psi$ be
 as in~(\ref{mappapsi}) and let's adopt the
 notation introduced in the proof of
 Lemma~(\ref{faine}). Let $\ell_0=\CC v_0\in
 Y^{sm}_A$. By Proposition~(\ref{yewgen}) we know that
 $\ell_0\notin W_A$ and hence
 $\ker\psi(\ell_0)=\CC v_0\wedge\alpha_0$ with
 $\alpha_0\in\wedge^2 V_0$ an indecomposable
 element; let $J_0\subset V_0$ be the span of
 $\alpha_0$, thus $\dim J_0=4$ because $\alpha_0$
  is indecomposable. Let $E_0\subset V$ be the
 codimension-$1$ subspace spanned by $v_0$ and
 $J_0$. It follows immediately
 from~(\ref{eccoderivata}) with $\beta_0=\alpha_0$ that
\begin{equation}\label{spaziotang}
  T_{\ell_0}Y_A=\PP(E_0).
\end{equation}
Now notice that
 $v_0\wedge\alpha_0\in\wedge^3 E_0\cap A$ and
 hence
\begin{equation}\label{voila}
 \{0\}\not=(\wedge^3 V/\wedge^3 E_0 + A)^{\vee}
=(\wedge^3 E_0)^{\bot}\cap A^{\bot}.
\end{equation}
Let $E_0^{\bot}=\CC \phi_0$; then
\begin{equation}\label{identifico}
  (\wedge^3 E_0)^{\bot}=\CC
\phi_0\otimes\wedge^2(V^{\vee}/\CC\phi_0)=
F_{\CC\phi_0}.
\end{equation}
By~(\ref{spaziotang})
and~(\ref{voila})-(\ref{identifico}) we get
that $T_{\ell_0}Y_A\in Y_{A^{\bot}}$; this
proves~(\ref{contenuto}).
 \end{proof}
By the above proposition duality defines a
rational involution on the set of projective
equivalence classes of EPW-sextics. We will
show later on - see Section~(\ref{involsuk}) -
that a generic EPW-sextic is not self-dual,
i.e.~the rational involution defined by duality
is not the identity.
 \section{Double covers of EPW-sextics}
 \label{epwdoppio}
 \setcounter{equation}{0}
We give the details of the following
observation: for $A\in \mathbb{LG}(\wedge^3
V)^0$ the variety $Y_A$ supports a quadratic
sheaf (see Definition~(0.2) of~\cite{cascat})
and if $X_A\to Y_A$ is the associated double
cover then $X_A$ is smooth. Let $A\in
\mathbb{LG}(\wedge^3 V)$ and let
$A^{\vee}\subset\wedge^3 V$ be a Lagrangian
subspace transversal to $A$ - see
Section~(\ref{eccoepw}). Thus we have
\begin{equation}\label{apiuaduale}
  \wedge^3 V=A\oplus A^{\vee}.
\end{equation}
Let
\begin{equation}
  \wt\lambda_A\colon
  \wedge^3 V\otimes\cO_{\PP(V)}\to
  A^{\vee}\otimes\cO_{\PP(V)},
  \end{equation}
be the projection corresponding to
Decomposition~(\ref{apiuaduale}). Let $\nu$ and
$\lambda_A$ be given by~(\ref{effeinwedge})
and~(\ref{mappaquot}) respectively; then
$\lambda_A=\wt\lambda_A\circ\nu$. We will study
the sheaf $coker(\lambda_A)$ fitting into the
exact sequence
\begin{equation}\label{sucker}
  0\to F\overset{\lambda_A}{\lra}
  A^{\vee}\otimes\cO_{\PP(V)}
  \lra coker(\lambda_A)\to 0.
\end{equation}
 \begin{prp}\label{presdiqua}
 Keep notation as above and assume that
 $A\in \mathbb{LG}(\wedge^3 V)^0$.
 Let $\ell_0\in \PP(V)$.
\begin{itemize}
  \item [(1)]
 If $\ell_0\notin Y_A$ then $coker(\lambda_A)$
 is zero in a neighborhood of $\ell_0$.
  \item [(2)]
 If $\ell_0\in (Y_A\setminus W_A)$
  there exist an open affine $U\subset\PP(V)$
 containing $\ell_0$ and trivializations
 $F|_U\cong\cO_U^9\oplus\cO_U$ and
 $A\otimes\cO_U\cong\cO_U^9\oplus\cO_U$ such
 that~(\ref{sucker}) restricted to $U$ reads
\begin{equation}\label{}
  0\lra \cO_U^9\oplus\cO_U
  \overset{(Id, x)}{\lra}\cO_U^9\oplus\cO_U\lra
  coker(\lambda_A)\otimes\cO_U\lra 0
\end{equation}
where $x$ is a local generator of $I_{Y_A\cap
U}$.
  \item [(3)]
If $\ell_0\in W_A$ there exist an open affine
$U\subset\PP(V)$
 containing $\ell_0$ and trivializations
 $F|_U\cong\cO_U^8\oplus\cO_U^2$ and
 $A\otimes\cO_U\cong\cO_U^8\oplus\cO_U^2$ such
 that~(\ref{sucker}) restricted to $U$ reads
\begin{equation}\label{(Id,M)}
  0\lra \cO_U^8\oplus\cO^2_U
  \overset{(Id,M)}{\lra}\cO_U^8\oplus\cO^2_U\lra
 coker(\lambda_A)\otimes\cO_U\lra 0
\end{equation}
where
\begin{equation}\label{matricem}
M=\left(\begin{array}{rr}
     x &  y \\
    y & z
  \end{array}\right),
\end{equation}
 with $x,y,z$ generators of $I_{W_A\cap U}$.
\end{itemize}
 \end{prp}
 \begin{proof}
 This is a straightforward consequence of the
 proof of Lemma~(\ref{faine}); we leave the
 details to the reader.
 \end{proof}
We will need a few results on the sheaf $coker(\lambda_A)$ and sheaves which
locally look like $coker(\lambda_A)$.
\begin{dfn}\label{fasciocc}
A  coherent sheaf $\cF$ on  a smooth projective variety $Z$ is a {\it
Casnati-Catanese sheaf} if for every $p\in Z$ there exists $U\subset Z$ open in
the classical topology containing $p$ such that one of the following holds:
\begin{itemize}
  \item [(1)]
 $\cF|_{U}=0$.
  \item [(2)]
 There exist $x\in Hol(U)$ with $x(p)=0$,
 $dx(p)\not=0$ and an exact sequence
\begin{equation}\label{primaris}
  0\to \cO_U\mapor{\cdot x}\cO_U\mapor{}
  \cF\to 0.
\end{equation}
  \item [(3)]
 There exist $x,y,z\in Hol(U)$ vanishing at $p$ with
 $dx(p),dy(p),dz(p)$ linearly independent and an
 exact sequence
\begin{equation}\label{secondaris}
  0\to \cO^2_U\mapor{M}\cO^2_U\mapor{}
  \cF\to 0
\end{equation}
  where $M$ is the map defined by
  Matrix~(\ref{matricem}).
\end{itemize}
\end{dfn}
Thus $coker(\lambda_A)$ is a typical example of
 a Casnati-Catanese sheaf. If $\cF$ is a
 Casnati-Catanese sheaf the schematic support of
$\cF$ is a divisor $D$ on $Z$; thus letting
$i\colon D\hra Z$ be the inclusion
we have
\begin{equation}\label{vienedad}
  \cF=i_{*}\cG
\end{equation}
for a coherent sheaf $\cG$ on $D$. In particular $coker(\lambda_A)=i_{*}\zeta_A$
for $i\colon Y_A\hra\PP(V)$ the inclusion map and
 $\zeta_A$  a certain coherent
sheaf on $Y_A$.
\begin{prp}\label{zetacona}
Let $\cF$ be a Casnati-Catanese sheaf on $Z$. Let $D$ be the schematic support
of $\cF$ and hence $\cF=i_{*}\cG$ where $i\colon D\hra Z$ is the inclusion,
see~(\ref{vienedad}).
\begin{itemize}
  \item [(1)]
If $q\ge 2$ then $Tor_q(\cF,\cE)=0$ for any
(abelian) sheaf $\cE$ on $Z$.
  \item [(2)]
The sheaf $\cG$  is locally isomorphic to
$\cG^{\vee}:=Hom(\cG,\cO_{D})$. In particular
$\cG$ is pure i.e.~there does not exist a
non-zero subsheaf of $\cG$ supported on a
subscheme of $D$ of dimension strictly smaller
than $\dim D$.
  \item [(3)]
 The map of sheaves $\cO_{D}\to Hom(\cG,\cG)$ which
 associates to
 $f\in\cO_{D,p}$ multiplication by $f$ is an isomorphism.
  \item [(4)]
 There is an isomorphism
\begin{equation}\label{dualeduale}
 Ext^1(\cF,\cO_Z)\cong
i_{*}\left(\cG^{\vee}\otimes N_{D/Z}\right).
\end{equation}
\end{itemize}
\end{prp}
\begin{proof}
(1) (2) and~(3) follow immediately from the
given local resolutions of $\cF$. (4): Let
$\cE_0\to\cF$ be a surjection with $\cE_0$
locally-free and let $\cE_1$ be the kernel of
the surjection. Thus we have an exact sequence
\begin{equation}\label{euno}
  0\to \cE_1\mapor{h}\cE_0\mapor{}\cF\to 0.
\end{equation}
By Item~(1) the sheaf $\cE_1$ is locally-free.
The dual of~(\ref{euno}) is the exact sequence
\begin{equation}\label{barbablu}
  0\to \cE_0^{\vee}\mapor{h}\cE_1^{\vee}
  \mapor{\partial}Ext^1(\cF,\cO_Z)\to 0.
\end{equation}
Multiplication defines an inclusion
\begin{equation}\label{primamappa}
  \cF\otimes\cO_Z(-D)=(\cE_0/\cE_1)\otimes\cO_Z(-D)
  \to\cE_0(-D)/\cE_1(-D).
\end{equation}
Since $\cF$ is supported on $D$ we have an
inclusion $\cE_0(-D)\hra\cE_1$ and hence an
inclusion
\begin{equation}\label{secondamappa}
  \cE_0(-D)/\cE_1(-D)\hra
  \cE_1/\cE_1(-D)=\cE_1\otimes\cO_D.
\end{equation}
Composing Map~(\ref{primamappa}) and
Map~(\ref{secondamappa}) we get an inclusion
\begin{equation}\label{}
  \cG\otimes\cO_D(-D)\hra\cE_1|_{D}
\end{equation}
whose dual is a surjection
\begin{equation}\label{cesiamo}
  \cE_1^{\vee}|_{D}\to \cG\otimes N_{D/Z}.
\end{equation}
Since $\cF$ is supported on $D$ the connecting
homomorphism map of~(\ref{barbablu})
annihilates $\cE^{\vee}_1(-D)$ and hence it may
be identified with a quotient map of
$\cE_1^{\vee}|_{D}$: a local computation shows
that the quotient map is~(\ref{cesiamo}). This
proves Item~(4).
\end{proof}
We set
\begin{equation}\label{defxia}
  \xi_A:=\zeta_A\otimes\cO_{Y_A}(-3).
\end{equation}
\begin{prp}\label{gianni}
Keep notation as above and assume that $A\in \mathbb{LG}(\wedge^3 V)^0$.
\begin{itemize}
  \item [(1)]
 There exists a symmetric isomorphism
\begin{equation}\label{mappazza}
  \alpha_A\colon\xi_A\overset{\sim}{\lra}\xi_A^{\vee}
\end{equation}
 defining a commutative multiplication map
\begin{equation}\label{prodottoxi}
  \overline\alpha_A\colon\xi_A\otimes\xi_A\lra\cO_{Y_A}.
\end{equation}
  \item [(2)]
Multiplication~(\ref{prodottoxi}) is an
isomorphism away from $W_A$ and near $\ell_0\in
W_A$ is described as follows. There exist an
open affine $U\subset Y_A$ containing $\ell_0$,
global generators $\{e_1,e_2\}$ of
$\xi_A\otimes\cO_U$ and $x,y,z\in\CC[U]$
generating the ideal of $W_A\cap U$   such that
\begin{equation}\label{formalocale}
  \overline\alpha_A(e_1\otimes e_1)=x,\quad
  \overline\alpha_A(e_1\otimes e_2)=
  \overline\alpha_A(e_2\otimes e_1)=y,\quad
  \overline\alpha_A(e_2\otimes e_2)=z.
\end{equation}
  \item [(3)]
Any map $\gamma\colon\xi_A\to\xi_A^{\vee}$ is a
 constant multiple of $\alpha_A$.
\end{itemize}
\end{prp}
\begin{proof}
Let
 \begin{equation}
 \wt\mu_A\colon
  \wedge^3 V\otimes\cO_{\PP(V)}\to
  A\otimes\cO_{\PP(V)}
 \end{equation}
be the projection given by
Decomposition~(\ref{apiuaduale}) and let
$\mu_A\colon F\to A\otimes\cO_{\PP(V)}$ be
defined by $\mu_A:=\wt\mu_A\circ\nu$ where
$\nu$ is as in~(\ref{effeinwedge}). The diagram
\begin{equation}
 \begin{array}{ccc}
F&\mapor{\lambda_A}&
A^{\vee}\otimes\cO_{\PP(V)}\\
\mapver{\mu_A}& &\mapver{\mu_A^{\vee}}\\
A\otimes\cO_{\PP(V)}& \mapor{\lambda_A^{\vee}}&
F^{\vee}
\end{array}
\end{equation}
is commutative because
$F\overset{(\mu_A,\lambda_A)}{\lra}(A\oplus
A^{\vee})\otimes\cO_{\PP(V)}$ is a Lagrangian
embedding. The map $\lambda_A$ is an injection
of sheaves because $Y_A\not=\PP(V)$ and hence
also $\lambda_A^{\vee}$ is an injection of
sheaves. Thus there is a unique $\beta_A\colon
i_{*}\zeta_A\lra
Ext^1(i_{*}\zeta_A,\cO_{\PP(V)})$ making the
following diagram commutative with exact
horizontal sequences:
\begin{equation}\label{spqr}
\begin{array}{ccccccccc}
0 & \to & F&\mapor{\lambda_A}& A^{\vee}\otimes\cO_{\PP(V)} & \lra & i_{*}\zeta_A
&
\to & 0\\
 & & \mapver{\mu_A}& &\mapver{\mu^{\vee}_A} &
&
\mapver{\beta_A}& & \\
0 & \to & A\otimes\cO_{\PP(V)}& \mapor{\lambda_A^{\vee}}& F^{\vee} & \lra &
Ext^1(i_{*}\zeta_A,\cO_{\PP(V)}) & \to & 0
\end{array}
\end{equation}
By Isomorphism~(\ref{dualeduale}) we get that
$\beta_A$ may be viewed as a map
$\beta_A\colon\zeta_A\overset{\sim}{\lra}\zeta^{\vee}_A(6)$.
\begin{clm}\label{essenziale}
Let $Z$ be a smooth projective variety and
$\cF$ be a Casnati-Catanese sheaf on $Z$.
Suppose that there exist vector-bundles
$\cE_0,\cE_1$ on $Z$ and an exact sequence
\begin{equation}\label{quadrato}
\begin{array}{ccccccccc}
0 & \to & \cE_1&\mapor{\lambda}& \cE_0 &
\mapor{} & \cF &
\to & 0\\
 & & \mapver{\mu}& &\mapver{\mu^{\vee}} &
&
\mapver{\beta}& & \\
0 & \to & \cE_0^{\vee}& \mapor{\lambda^{\vee}}&
\cE_1^{\vee} & \mapor{} & Ext^1(\cF,\cO_Z) &
\to & 0
\end{array}
\end{equation}
Then $\beta$ is an isomorphism if and only if
the map
\begin{equation}\label{}
  \cE_1\mapor{(\lambda,\mu)}\cE_0\oplus\cE_0^{\vee}
\end{equation}
is an injection of vector-bundles i.e.~it is
injective on fibers.
\end{clm}
\begin{proof}
Let $\cF=i_{*}\cG$ where $i\colon D\hra Z$ is the inclusion. First we notice
that $\beta$ is an isomorphism if and only if it is surjective; in fact by
Items~(4) and~(2) of Proposition~(\ref{zetacona}) we have local identifications
of the sheaves $Hom(\cF,Ext^1(\cF,\cO_Z))$ with $Hom(\cG,\cG)$ and it follows
from Item~(3) of the same proposition that a map of stalks $\cG_p\to\cG_p$ is an
isomorphism if and only if it is surjective. By Nakayama's Lemma (or by a direct
computation) we get that $\beta$ is an isomorphism if and only if for every
$p\in Z$ the map from the fiber of $\cF$ at $p$ to the fiber of
$Ext^1(\cF,\cO_Z)$ at $p$ is surjective. As is easily checked there exist
$U\subset Z$ open in the classical topology containing $p$, trivial
vector-bundles $\cA_i,\cB_i$ on $U$ for $i=0,1$ and isomorphisms $\cE_i|_U\cong
\cA_i\oplus\cB_i$ such that the restriction of~(\ref{quadrato}) to $U$ reads
\begin{equation}\label{}
\begin{array}{ccccccccc}
0 & \to & \cA_1\oplus\cB_1&\mapor{(\phi,\psi)}& \cA_0\oplus\cB_0 & \mapor{} &
\cF|_U &
\to & 0\\
 & & \mapver{\mu}& &\mapver{\mu^{\vee}} &
&
\mapver{\beta}& & \\
0 & \to & \cA_0^{\vee}\oplus\cB_0^{\vee}& \mapor{(\phi^{\vee},\psi^{\vee})}&
 \cA_1^{\vee}\oplus\cB_1^{\vee} &
\mapor{} & Ext^1(\cF_U,\cO_U) & \to & 0
\end{array}
\end{equation}
where $\phi\colon\cA_1\to\cA_0$ is an isomorphism and $\psi\colon\cB_1\to\cB_0$
is a standard resolution of $\cF_U$ as given in Definition~(\ref{fasciocc}),
i.e.~$\cB_1=\cB_0=0$ if~(1) of Definition~(\ref{fasciocc}) holds,
$\cB_1\cong\cB_0\cong\cO_U$  and $\psi$ is the map of~(\ref{primaris}) if~(2) of
Definition~(\ref{fasciocc}) holds and finally if~(3) of
Definition~(\ref{fasciocc}) holds then $\cB_1\cong\cB_0\cong\cO^2_U$ and $\psi$
is the map of~(\ref{secondaris}). The map $\beta$ is a surjection from the fiber
of $\cF$ at $p$ to the fiber of $Ext^1(\cF_U,\cO_U)$ if and only if the
composition
\begin{equation}\label{primacomp}
  \cB_0\mapor{\iota}\cA_0\oplus\cB_0\mapor{\mu^{\vee}}
  \cA^{\vee}_1\oplus\cB^{\vee}_1
 \mapor{\pi}\cB^{\vee}_1
\end{equation}
is surjective at $p$. On the other hand since $\psi$ is zero  at $p$ the map
$((\phi,\psi),\mu)$  is injective  at $p$ if and only if the composition
\begin{equation}\label{secondacomp}
  \cB_1\mapor{\pi^{\vee}}\cA_1\oplus\cB_1\mapor{\mu}
  \cA^{\vee}_0\oplus\cB^{\vee}_0
\end{equation}
is injective at $p$. Since $\psi$ vanishes at $p$ while $\phi^{\vee}$ is an
isomorphism at $p$ the equality
$(\phi^{\vee},\psi^{\vee})\circ\mu=\mu^{\vee}\circ(\phi,\psi)$ gives that the
composition
\begin{equation}\label{}
 \cB_1\mapor{\pi^{\vee}}\cA_1\oplus\cB_1\mapor{\mu}
  \cA^{\vee}_0\oplus\cB^{\vee}_0 \mapor{}\cA^{\vee}_0
\end{equation}
vanishes at $p$. Thus~(\ref{secondacomp}) is injective at $p$ if and only if the
composition
\begin{equation}\label{terzacomp}
  \cB_1\mapor{\pi^{\vee}}\cA_1\oplus\cB_1\mapor{\mu}
  \cA^{\vee}_0\oplus\cB^{\vee}_0\mapor{\iota^{\vee}}\cB^{\vee}_0
\end{equation}
is injective at $p$. Since Composition~(\ref{terzacomp}) is the transpose of
 Composition~(\ref{primacomp}) this proves that $\beta$ is a surjection at $p$
 if and only if $((\phi,\psi),\mu)$  is injective  at $p$.
\end{proof}
The above result together with Exact
Sequence~(\ref{spqr}) gives that $\beta_A$ is
an isomorphism. We define the map $\alpha_A$
of~(\ref{mappazza}) to be the tensor product of
$\beta_A$ times the identity map on
$\cO_{Y_A}(-3)$. Since $\beta_A$  is an
isomorphism we get that $\alpha_A$ is an
isomorphism . The map
$\alpha_A^{\vee}-\alpha_A$ is zero on
$(Y_A\setminus W_A)$ because $\xi_A$ is
locally-free of rank $1$ on $(Y_A\setminus
W_A)$. By Item~(2) of
Proposition~(\ref{zetacona}) $\zeta_A$ is pure
and hence $\xi_A$ is pure; thus
$\alpha_A^{\vee}-\alpha_A$ is zero on $Y_A$
i.e.~$\alpha_A$ is symmetric. This proves
Item~(1). Let's prove Item~(2). Since $\xi_A$
is locally-free of rank $1$ away from $W_A$ it
is clear that~(\ref{prodottoxi}) is an
isomorphism away from $W_A$.
Equation~(\ref{formalocale}) defines a
symmetric isomorphism
$\mu\colon\xi_A\otimes\cO_U\to
\xi^{\vee}_A\otimes\cO_U$. By Item~(3) of
Proposition~(\ref{zetacona}) any other
symmetric isomorphism $\xi_A\otimes\cO_U\to
\xi^{\vee}_A\otimes\cO_U$ is equal to
$f\cdot\mu$ where $f\in\CC[U]^{\times}$ is
invertible; Item~(2) of the proposition follows
at once from this.
 Lastly we prove Item~(3).
 The composition $\alpha_A^{-1}\circ\gamma$ is a
 global section
 of $Hom(\xi_A,\xi_A)$ and by Item~(3) of
 Proposition~(\ref{zetacona}) we get that
$\alpha_A^{-1}\circ\gamma$ is multiplication by
a constant $c$; thus $\gamma=c\alpha_A$.
\end{proof}
Now we are ready to define the double cover of
$X_A$ for $A\in \mathbb{LG}(\wedge^3 V)^0$.
Map~(\ref{prodottoxi}) gives
$\cO_{Y_A}\oplus\xi_A$ the structure of a
commutative finite $\cO_{Y_A}$-algebra. Let
\begin{equation}\label{doppioepw}
  X_A:=Spec\left(\cO_{Y_A}\oplus\xi_A\right)
\end{equation}
and $f\colon X_A\to Y_A$ be the structure map.
Thus $f$ is a finite map of degree $2$; let
$\phi\colon X_A\to X_A$ be the covering
involution - thus $\phi$ corresponds to the map
on $\cO_{Y_A}\oplus\xi_A$ which is the identity
on $\cO_{Y_A}$ and multiplication by $(-1)$ on
$\xi_A$.  Equivalently
\begin{equation}\label{}
  f_{*}\cO_{X_A}=\cO_{Y_A}\oplus\xi_A
\end{equation}
with $\cO_{Y_A}$ and $\xi_A$ the $+1$ and
$(-1)$-eigensheaf respectively. It follows at
once from our lemma that $f$ is ramified
exactly over $W_A$ and that $X_A$ is smooth.
Let $A$ vary in $\mathbb{LG}(\wedge^3 V)^0$. By
Item~(2) of Remark~(\ref{modelloloc}) the
family of $Y_A$'s  is locally trivial; since
$\mathbb{LG}(\wedge^3 V)^0$ is irreducible we
get the following result.
\begin{prp}\label{tuttedef}
The varieties $X_A$ for $A$ varying
in $\mathbb{LG}(\wedge^3 V)^0$ are all
deformation equivalent.
\end{prp}
 \section{Proof of Theorem~(\ref{mainthm})}
 \label{dimprinc}
 \setcounter{equation}{0}
Throughout this section $X$ is a numerical
$(K3)^{[2]}$ with an ample divisor $H$ such that
both~(\ref{gradodod}) and Item~(a) of
Section~(\ref{rpologo}) hold. Thus we have  an
anti-symplectic involution $\phi\colon X\to X$
with quotinet map $f\colon X\to
X/\langle\phi\rangle=:Y$ and an embedding
\begin{equation}\label{ysestica}
  j\colon Y\hra |H|^{\vee},\quad
  \deg Y=6
\end{equation}
such that $j\circ f$ is the tautological map
$X\to |H|^{\vee}$. Let $X^{\phi}$  be the fixed
locus of $\phi$; then $sing Y\cong X^{\phi}$.
Since $\phi$ is anti-symplectic $X^{\phi}$ is a
smooth Lagrangian surface - not empty because a
smooth hypersurface in $\PP^5$ is simply
connected. Thus
\begin{equation}\label{pusingy}
  \dim (sing Y)=2
\end{equation}
and $f$ is unramified over $(Y\setminus sing Y)$. We have a decomposition
\begin{equation}\label{eigendec}
  f_{*}\cO_X=\cO_Y\oplus\eta,
\end{equation}
where $\eta$ is the $(-1)$-eigensheaf of the involution $\phi^{*}\colon
f_{*}\cO_X\to f_{*}\cO_X$.
\begin{clm}\label{itiscc}
Keep notation as above. The sheaf $j_{*}\eta$ is a Casnati-Catanese sheaf on
$|H|^{\vee}$.
\end{clm}
\begin{proof}
Let $p\in (|H|^{\vee}\setminus Y)$; then $j_{*}\eta$ is zero in a neighborhood
of $p$. Let $p\in(Y\setminus sing Y)$; then $\eta$ is locally-free of rank $1$
in a neighborhood of $p$ (in $Y$) and hence we see that~(2) of
Definition~(\ref{fasciocc}) holds. Let $p\in sing Y$. Let
$f^{-1}(p)=\{\wt{p}\}$. Since $\dim X^{\phi}=2$ there exist a $\wt{U}\subset X$
open in the classical topology containing $\wt{p}$ with analytic coordinates
$\{u_1,u_2,v_1,v_2\}$ centered at $\wt{p}$  and an open $U\subset |H|^{\vee}$
containing $p$ with analytic coordinates $\{x_1,\ldots,x_5\}$ centered at $p$
such that
\begin{align}
\phi^{*}(u_1,u_2,v_1,v_2) & =  (-u_1,-u_2,v_1,v_2)\label{mappaphi} \\
  f^{*}(x_1,x_2,x_3,x_4,x_5)& =
  (u^2_1,u^2_2,-u_1 u_2,v_1,v_2).\label{}
\end{align}
It follows from~(\ref{mappaphi}) that
\begin{equation}\label{generoj}
  \text{$j_{*}\eta|_U$ is generated by
  $(u_1,u_2)$.}
\end{equation}
A straightforward computation gives the free
presentation
\begin{equation}\label{marcodemarchi}
  0\to \cO_U^2
  \mapor{\left(\begin{array}{rr}
     x_2 &  x_3 \\
    x_3 & x_1
  \end{array}\right)}
  \cO_U^2
  \mapor{\left(u_2,u_1\right)}
  (j_{*}\eta)|_U\to 0.
\end{equation}
Thus~(3) of Definition~(\ref{fasciocc}) holds.
\end{proof}
For future use we notice the following: keeping notation as in the above proof
 we may assume by shrinking $U$ that $f^{-1}(U)=\wt{U}$ and then
\begin{equation}\label{eqloc}
  Y\cap U=V(x_1 x_2-x_3^3).
\end{equation}
Multiplication on $f_{*}\cO_X$ defines a symmetric map
\begin{equation}\label{moltsueta}
  \ov{\alpha}\colon\eta\otimes\eta\to\cO_Y
\end{equation}
making $\eta$ a quadratic sheaf in the sense of
Casnati-Catanese~\cite{cascat}. A
straightforward computation shows that
$\ov{\alpha}$ defines a symmetric isomorphism
\begin{equation}\label{isoalfa}
  \alpha\colon\eta\overset{\sim}{\lra}
  \eta^{\vee}:=Hom(\eta,\cO_Y).
\end{equation}
The symmetric map~(\ref{moltsueta}) makes
$\cO_Y\oplus\eta$ a commutative $\cO_Y$-algebra
and we have the tautological isomorphism
\begin{equation}\label{ixspec}
  X\cong Spec(\cO_Y\oplus\eta).
\end{equation}
The main result of this section is the
following.
\begin{thm}\label{chiappone}
Keep notation and hypotheses as above. There exists $A\in \mathbb{LG}(\wedge^3
V)^0$ such that $Y=Y_A$ and $\eta\cong\xi_A$.
\end{thm}
The proof of Theorem~(\ref{chiappone}) will be
given at the end of this section. Let's show that
Theorem~(\ref{mainthm}) follows  from
Theorem~(\ref{chiappone}). Let's prove~(1) of
Theorem~(\ref{mainthm}) i.e.~that $f\colon X\to
Y$ is identified with the natural double cover
$X_A\to Y_A$. According to
Theorem~(\ref{chiappone}) we may identify $\eta$
and $\xi_A$. We claim that with this
identification the map $\alpha$
of~(\ref{isoalfa}) gets identified with a
non-zero constant multiple of $\alpha_A$; in fact
$\alpha$ is non-zero and hence the claim follows
from Item~(3) of Proposition~(\ref{gianni}). Thus
multiplying the isomorphism
$\eta\overset{\sim}{\lra}\xi_A$ by a suitable
constant we may assume that the map $\alpha$ gets
identified with $\alpha_A$; by~(\ref{ixspec})
and~(\ref{doppioepw}) we get that $f\colon X\to
Y$ is identified with $X_A\to Y_A$. Let's
prove~(2) of Theorem~(\ref{chiappone}) i.e.~that
if $A\in\mathbb{LG}(\wedge^3 V)^0$ then $X_A\to
Y_A$ is a deformation of $(K3)^{[2]}$. By
Proposition~(\ref{tuttedef}) the $X_A$ for
$A\in\mathbb{LG}(\wedge^3 V)^0$ are all
deformation equivalent. Since every K\"ahler
deformation of an irreducible symplectic manifold
is an irreducible symplectic manifold it suffices
to prove that there exists one
$A\in\mathbb{LG}(\wedge^3 V)^0$ such that $X_A$
which is a symplectic irreducible variety
deformation equivalent of $(K3)^{[2]}$. By
Item~(1) of Theorem~(\ref{mainthm}) it suffices
to exhibit  $X,H$ where $X$ is a deformation of
$(K3)^{[2]}$ and $H$ is an ample divisor on $X$
such that both~(\ref{gradodod}) and Item~(a) of
Section~(\ref{prologo}) hold. Such an example was
given by Mukai (Ex.(5.17) of~\cite{mukai}),
details are in Subsection~(5.4) of~\cite{oginv}.
We briefly describe the example. Let
$F\subset\PP^6$ be \lq\lq the\rq\rq Fano $3$-fold
of index $2$ and degree $5$ i.e.~the transversal
intersection of $Gr(2,\CC^5)\subset\PP(\wedge^2
\CC^5)$ and a $6$-dimensional linear subspace of
$\PP(\wedge^2\CC^5)$. Let $\ov{Q}\subset\PP^6$ be
a quadric hypersurface intersecting transversely
$F$ and let
\begin{equation}\label{kappatre}
  S:=F\cap\ov{Q}.
\end{equation}
 Thus $(S,\cO_S(1))$ is a
generic polarized $K3$ surface of degree $10$. We
assume that
\begin{equation}\label{dieci}
  E\cdot c_1(\cO_S(1))\equiv 0\pmod{10},\quad
  \text{$E$ divisor on $S$,}
\end{equation}
this is Hypothesis~(4.8) of~\cite{oginv} and it
holds for $\ov{Q}$ belonging to an open dense
subset of the space of all quadrics. We have
$\dim|I_F(2)|=4$ and $\dim|I_S(2)|=5$. Let
$\Sigma$ be the degree-$7$ divisor on $|I_S(2)|$
parametrizing singular quadrics. Every
$Q\in|I_F(2)|$ is singular and hence
\begin{equation}\label{eccoy}
  \Sigma=|I_S(2)|+Y.
\end{equation}
Since $\deg\Sigma=7$ we have $\deg Y=6$. The
generic $Q\in Y$ has corank $1$ and hence it
has two rulings by $3$-dimensional linear
spaces. Thus there exists a natural double
cover of an open dense subset of $Y$. It turns
out (see Ex.(5.17) of~\cite{mukai}) and
Subsection~(5.4) of~\cite{oginv}) that the
double cover extends to a double cover $f\colon
X\to Y$, that
\begin{equation}\label{eccox}
  X=\{\text{$F$ stable sheaf on $S$,
   $rk(F)=2$, $c_1(F)=c_1(\cO_S(1)$, $c_2(F)=5$}\}/
   \text{isomorphism}
\end{equation}
and that $X$ is a  deformation of $S^{[2]}$. Let
$H:=f^{*}\cO_Y(1)$. As shown in Subsection~(5.4)
of~\cite{oginv} both~(\ref{gradodod}) and
Item~(a) of Section~(\ref{prologo}) hold. This
proves that~(2) of Theorem~(\ref{mainthm}) holds.
\subsection{A locally free resolution of
$j_{*}(\eta\otimes\cO_Y(3))$}
\begin{prp}\label{moltoampio}
Let $X$ be a numerical $(K3)^{[2]}$ with an ample
divisor $H$ such that both~(\ref{gradodod}) and
Item~(a) of Section~(\ref{prologo}) hold. If
$k\ge 3$ then $\cO_X(kH)$ is very ample.
\end{prp}
\begin{proof}
By hypothesis $nH$ is very ample for some $n\gg
0$ and hence it suffices to prove that the
multiplication map
\begin{equation}\label{moltsuix}
  H^0(\cO_X(kH))\otimes H^0(\cO_X(H))\lra
  H^0(\cO_X((k+1)H))
\end{equation}
is surjective for $k\ge 3$. By~(\ref{pusingy})
there exists a plane $\Lambda\subset
|H|^{\vee}$ such that $\Lambda\cap (sing
Y)=\emptyset$ and $\Lambda$ is transversal to
$Y$. Let $C:=\Lambda\cap Y$ and
$\wt{C}:=f^{-1}C$. Then $C$ is a smooth plane
sextic by~(\ref{ysestica}) and $\pi:=f|_C\colon
\wt{C}\to C$ is an unramified double cover.
Let's prove that
\begin{equation}\label{moltsuct}
  H^0(\cO_{\wt{C}}(kH))\otimes H^0(\cO_{\wt{C}}(H))\lra
  H^0(\cO_{\wt{C}}((k+1)H))
\end{equation}
is surjective for $k\ge 3$. We have
\begin{equation}\label{}
  \pi_{*}\cO_{\wt{C}}=\cO_C\oplus \lambda
\end{equation}
where $\lambda$ is a non-trivial square root of
$\cO_{\wt{C}}$. Thus
\begin{equation}\label{}
  H^0(\cO_{\wt{C}}(kH))=H^0(\pi_{*}\cO_{\wt{C}})=
  H^0(\cO_C(k))\oplus H^0(\lambda(k)).
\end{equation}
Since the multiplication map
$H^0(\cO_C(k))\otimes H^0(\cO_C(1))\to
H^0(\cO_C(k+1))$ is surjective it suffices to
prove surjectivity of
\begin{equation}\label{bastacio}
H^0(\lambda(k))\otimes H^0(\cO_C(1))\to
H^0(\lambda(k+1))
\end{equation}
for $k\ge 3$. Since $C$ is a smooth plane
sextic adjunction gives $K_C\cong\cO_C(3)$;
since $\lambda\not\cong\cO_C$ we get that
$h^1(\lambda(k))=0$ for $k\ge 3$. Thus
\begin{equation}\label{eccodim}
  h^0(\lambda(k))=\chi(\lambda(k))=6k-9,\quad k\ge 3.
\end{equation}
Let $U\subset H^0(\cO_C(1))$ be a
$2$-dimensional subspace spanned by sections
$\epsilon_0,\epsilon_1$ with no common zeroes.
Consider the multiplication map
\begin{equation}\label{uperacca}
  H^0(\lambda(k))\otimes U
  \overset{\mu}{\lra} H^0(\lambda(k+1)).
\end{equation}
By the base-point-free pencil trick we have
\begin{equation}\label{}
  \ker(\mu)=\{(\sigma\epsilon_0)\otimes \epsilon_1-
  (\sigma\epsilon_1)\otimes \epsilon_0|\
  \sigma\in H^0(\lambda(k-1))\}.
\end{equation}
Using Formula~(\ref{eccodim}) one gets that
$\dim Im(\mu)=h^0(\lambda(k+1))$ and hence
$\mu$ is surjective: thus Map~(\ref{bastacio})
is surjective and hence also
Map~(\ref{moltsuct}) is surjective. Now we
prove that Map~(\ref{moltsuix}) is surjective.
Let $X=X_4\supset X_3\supset X_2\supset
X_1=\wt{C}$ be a chain of smooth linear
sections of $X$, i.e.~$X_3\in |H|$ and
$X_2=D\cap D'$ where $D,D'\in |H|$ intersect
transversely. We claim that the restriction map
\begin{equation}\label{pappa}
  H^0(\cO_{X_i}(sH))\to H^0(\cO_{X_{i-1}}(sH))
\end{equation}
is surjective for $s=1$ and $s\ge 3$. It
suffices to show that
\begin{equation}\label{}
  h^1(\cO_{X_i}((s-1)H)),\quad s=1,\quad s\ge
  3.
\end{equation}
This follows from the Lefschetz Hyperplane
Section Theorem and Kodaira Vanishing.
Surjectivity of Map~(\ref{pappa}) for $s=1$ and
$s\ge 3$ together with surjectivity
of~(\ref{moltsuct}) gives surjectivity
of~(\ref{moltsuix}) by an easy well-known
argument.
\end{proof}
Let
\begin{equation}\label{tetadef}
  \theta:=\eta\otimes\cO_Y(3).
\end{equation}
\begin{crl}\label{treacca}
  Keep notation and hypotheses as above. Then
  $\theta$ is globally generated.
\end{crl}
\begin{proof}
Let $H^0(\cO_X(3H))^{-}\subset H^0(\cO_X(3H))$
be the (-1)-eigenspace for the action of
$\phi^{*}$. Then
\begin{equation}\label{tretregiugiu}
  H^0(\cO_X(3H))=H^0(\cO_Y(3H))\oplus
  H^0(\theta),\quad
  H^0(\theta)=H^0(\cO_X(3H))^{-}.
\end{equation}
Let $p\in Y$ and let $\theta_p$ be the fiber of
$\theta$ at $p$; we must show that
\begin{equation}\label{valuto}
  H^0(\cO_X(3H))^{-}\lra \theta_p
\end{equation}
is surjective. Assume that $p\in(Y\setminus
sing Y)$. Let $f^{-1}(p)=\{p_1,p_2\}$ and let
$L_{p_i}$ be the fiber of $\cO_X(3H)$ at $p_i$.
Since $\cO_X(3H)$ is very ample the evaluation
map $H^0(\cO_X(3H))\to(L_{p_1}\oplus L_{p_2})$
is surjective. Since $\theta_p$ is identified
with the $(-1)$-eigenspace for the action of
$\phi^{*}$ on $(L_{p_1}\oplus L_{p_2})$ we get
that~(\ref{valuto}) is surjective. Finally
assume that $p\in sing Y$ and let
$f^{-1}(p)=\{\wt{p}\}$. Then~(\ref{generoj})
identifies $\theta_p$ with
$\Omega^1_{X,\wt{p}}$ and Map~(\ref{valuto})
with
 differentiation at $\wt{p}$. Since
$\cO_X(3H)$ is very ample the differential at
$\wt{p}$ of the map $X\to |3H|^{\vee}$ is
injective; it follows that~(\ref{valuto}) is a
surjection.
\end{proof}
Let $\epsilon\colon
H^0(\theta)\otimes\cO_{|H|^{\vee}}\to
j_{*}\theta$ be the evaluation map. By the
above corollary $\epsilon$ is surjective; let
$G$ be the kernel of $\epsilon$. Thus we have
an exact sequence
\begin{equation}\label{succdigi}
  0\to G\lra H^0(\theta)\otimes\cO_{|H|^{\vee}}
  \overset{\epsilon}{\lra} j_{*}\theta\to 0.
\end{equation}
\begin{prp}\label{difftre}
Keep notation and assumptions as above. Then
\begin{equation}\label{}
  G\cong \Omega^3_{|H|^{\vee}}(3).
\end{equation}
\end{prp}
\begin{proof}
First we prove that $G$ is locally-free. By Claim~(\ref{itiscc}) the sheaf
$j_{*}\eta$ is Casnati-Catanese and hence so is $j_{*}\theta$. By~(1) of
Proposition~(\ref{zetacona}) we get that $G$ is locally-free.
 Since $G$ is
locally-free Beilinson's spectral sequence with
\begin{equation}\label{}
  E_1^{p,q}=H^q(G(p))\otimes
  \Omega_{|H|^{\vee}}^{-p}(-p)
\end{equation}
converges in degree $0$ to the graded sheaf
associated to a filtration on $G$,
see~\cite{oss} p.~240. Thus it suffices to
prove that
\begin{equation}\label{accagipi}
h^q(G(p))=
  \begin{cases}
    0 & \text{if $-5\le p\le 0$ and
    $(p,q)\not=(-3,3)$}, \\
    1 & \text{if $(p,q)=(-3,3)$}.
  \end{cases}
\end{equation}
This follows from a straighforward computation
which goes as follows.
Tensorizing~(\ref{succdigi}) by
$\cO_{|H|^{\vee}}(p)$ and taking the associated
cohomology exact sequence we get
\begin{multline}\label{succlunga}
  \cdots \to H^0(\theta)\otimes
  H^{q-1}(\cO_{|H|^{\vee}}(p))\lra
  H^{q-1}(\theta(p))\overset{\partial}{\lra} \\
  \overset{\partial}{\lra} H^q(G(p))\to
 H^0(\theta)\otimes
  H^{q}(\cO_{|H|^{\vee}}(p))\lra\cdots
\end{multline}
(We let $\theta(p):=\theta\otimes\cO_Y(p)$.)
From this one easily gets that
\begin{equation}\label{accaunogi}
 h^1(G)=0, \quad
 h^0(G(p)), \quad -5\le p\le 0,
\end{equation}
and that
\begin{equation}\label{deliso}
 h^{q-1}(\theta(p))=
 h^q(G(p)), \quad -5\le p\le 0,\
 (p,q)\not=(0,1).
\end{equation}
We compute the left-hand side. The map $f$ is
finite and we
have~(\ref{eigendec})-~(\ref{tetadef}); thus
\begin{multline}\label{findecomp}
 h^{q-1}(\cO_X((3+p)H))=h^{q-1}(f_{*}\cO_X((3+p)H))=\\
 =h^{q-1}(\cO_Y(3+p))\oplus
 h^{q-1}(\theta(p)).
\end{multline}
In order to compute $h^{q-1}(\theta(p))$ we
first compute $h^{q-1}(\cO_X((3+p)H))$. Kodaira
Vanishing gives that $h^{q-1}(\cO_X((3+p)H))=0$
for $2\le q\le 5$ and $p\not=-3$ and hence
by~(\ref{findecomp}) we get that
\begin{equation}\label{primaform}
 h^{q-1}(\theta(p))=0,\quad 2\le q\le 5,\
 p\not=-3.
\end{equation}
Now consider $q=1$. We claim that
\begin{equation}\label{tuttopari}
  h^0(\cO_X((3+p)H))=h^0(\cO_Y((3+p)H)),\quad
  -5\le p\le -1.
\end{equation}
For $-5\le p\le -3$ the equation is trivial and for $p=-2$ it holds by
hypothesis. To check equality for $p=-1$ we apply Formula~(4.0.4)
of~\cite{ognum}:
\begin{equation}
\chi(\cO_X(nH))=
\frac{1}{2}n^4+\frac{5}{2}n^2+3.
\end{equation}
Since $H$ is ample Kodaira Vanishing gives that
$h^0(\cO_X(2H))=\chi(\cO_X(2H))$ and by the
above formula we get $h^0(\cO_X(2H))=21$. On
the other hand a straightforward computation
gives that $h^0(\cO_Y(2))=21$ and
by~(\ref{findecomp}) this finishes the proof
of~(\ref{tuttopari}). Thus we have proved that
\begin{equation}\label{secondaform}
 h^0(\theta(p))=0,\quad -5\le p\le -1.
\end{equation}
Finally consider $h^{q-1}(\cO_X)$: it vanishes
for $q=2,4$ and
\begin{equation}\label{}
  h^{q-1}(\cO_X)=1=h^{q-1}(\cO_Y),\quad q=1,5.
\end{equation}
On the other hand $h^2(\cO_X)=1$ and
$h^2(\cO_Y)=0$. This proves that
\begin{equation}\label{terzaform}
  h^{q-1}(\theta(-3))=
  \begin{cases}
    0 & \text{if $q\not=3$}, \\
    1 & \text{if $q=3$}.
  \end{cases}
\end{equation}
Equation~(\ref{accagipi}) follows from
Formulae~(\ref{accaunogi}), (\ref{deliso}),
(\ref{primaform}), (\ref{secondaform})
and~(\ref{terzaform}).
\end{proof}
By the above proposition and
by~(\ref{succdigi}) we have an exact sequence
\begin{equation}\label{kappasucc}
  0\to \Omega^3_{|H|^{\vee}}(3)\overset{\kappa}\lra
  H^0(\theta)\otimes\cO_{|H|^{\vee}}
  \overset{\epsilon}{\lra} j_{*}\theta\to 0.
\end{equation}
\subsection{Proof of Theorem~(\ref{chiappone})}
\begin{clm}\label{isoditeta}
Keep notation as above. There exists an
isomorphism
\begin{equation}\label{isoteta}
  \beta\colon
  j_{*}\theta\overset{\sim}{\lra}
  Ext^1
  (j_{*}\theta,\cO_{|H|^{\vee}}).
\end{equation}
\end{clm}
\begin{proof}
 Since $j_{*}\theta$ is a Casnati-Catanese sheaf  we have an
 isomorphism
\begin{equation}\label{}
  Ext^1
  (j_{*}\theta,\cO_{|H|^{\vee}})\cong
 j_{*}\left(\theta^{\vee}
 \otimes N_{Y/|H|^{\vee}}\right)
\end{equation}
because of Item~(4) of Proposition~(\ref{zetacona}). By~(\ref{tetadef}) and
Isomorphism~(\ref{isoalfa}) we get
\begin{equation}\label{}
  j_{*}\left(\theta^{\vee}
 \otimes N_{Y/|H|^{\vee}}\right)=
 j_{*}\left(\eta^{\vee}\otimes\cO_Y(-3)
 \otimes N_{Y/|H|^{\vee}}\right)\cong
 j_{*}\left(\eta\otimes\cO_Y(-3)
 \otimes N_{Y/|H|^{\vee}}\right).
\end{equation}
By~(\ref{ysestica}) we have
$N_{Y/|H|^{\vee}}\cong\cO_Y(6)$ and hence we
get that
\begin{equation}\label{}
  j_{*}\left(\eta\otimes\cO_Y(-3)
 \otimes N_{Y/|H|^{\vee}}\right)\cong
 j_{*}\left(\eta\otimes\cO_Y(3)\right)=j_{*}\theta.
\end{equation}
The above equations prove~(\ref{isoteta}).
\end{proof}
Let $\kappa,\epsilon$ be as
in~(\ref{kappasucc}) and $\beta$ be as
in~(\ref{isoteta}). We claim that there exists
a map
\begin{equation}\label{}
  s\colon
H^0(\theta)\otimes\cO_{|H|^{\vee}}\to
\left(\Omega^3_{|H|^{\vee}}(3)\right)^{\vee}=
\Theta^3_{|H|^{\vee}}(-3)
\end{equation}
such that the following diagram is commutative:
\begin{equation}\label{grancomm}
 \begin{array}{ccccccccc}
0 & \mapor{} & \Omega^3_{|H|^{\vee}}(3)
&\mapor{\kappa}&
H^0(\theta)\otimes\cO_{|H|^{\vee}}  &
\mapor{\epsilon} & j_{*}\theta &
\to & 0\\
 & & \mapver{s^{\vee}}& &\mapver{s} &
&
\mapver{\beta}& & \\
0 & \mapor{} &
H^0(\theta)^{\vee}\otimes\cO_{|H|^{\vee}} &
\mapor{\kappa^{\vee}}&
\Theta^3_{|H|^{\vee}}(-3) & \mapor{\partial} &
Ext^1(j_{*}\theta,\cO_{|H|^{\vee}}) & \mapor{}
& 0
\end{array}
\end{equation}
In fact this follows from the results of
Casnati-Catanese~\cite{cascat} or of
Eisenbud-Popescu-Walter~\cite{epw}: by the
proof of Claim~(2.1) of~\cite{cascat} the
obstruction to existence of $s$ lies in
$H^1(Sym_2
\left(H^0(\theta)\otimes\cO_{|H|^{\vee}}
\right)^{\vee})$ which is zero and hence $s$
exists.
\begin{rmk}
{\rm Proposition~(1.6) of~\cite{cascat} does
not hold with $\cF=j_{*}\theta$ because
$\chi(j_{*}\theta(-3))$ is not even, see
Theorem~(9.1) of~\cite{epw} - in fact
$\chi(j_{*}\theta(-3))=1$. Thus unlike the
surfaces considered by Casnati-Catanese the
$4$-fold $Y$ cannot be presented as the
degeneracy locus of a symmetric map of
vector-bundles.}
\end{rmk}
\begin{clm}\label{iniezione}
Keep notation and assumptions as above. Then
\begin{equation}\label{}
 \Omega^3_{|H|^{\vee}}(3)\mapor{(\kappa,s^{\vee})}
 \left(H^0(\theta)\oplus
 H^0(\theta)^{\vee}\right)
 \otimes\cO_{|H|^{\vee}}
\end{equation}
is an injection of vector-bundles. The image of $(\kappa,s^{\vee})$ is
Lagrangian for the tautological symplectic form  on $H^0(\theta)\oplus
 H^0(\theta)^{\vee}$ given by
\begin{equation}\label{formataut}
  \lambda((\alpha,\psi),(\alpha',\psi')):=
  \psi(\alpha)-\psi'(\alpha').
\end{equation}
\end{clm}
\begin{proof}
The sheaf $j_{*}\theta$ is a Casnati-Catanese
sheaf on $|H|^{\vee}$; since $\beta$ is an
isomorphism we get by Claim~(\ref{essenziale})
that $(\kappa,s^{\vee})$ is an injection of
vector-bundles.   The tautological symplectic
form vanishes on $Im(\kappa,s^{\vee})$ by
commutativity of Diagram~(\ref{grancomm}).
Since $\Omega^3_{|H|^{\vee}}(3)$ has rank $10$
it follows that $Im(\kappa,s^{\vee})$ is
Lagrangian.
\end{proof}
We will show that Diagram~(\ref{grancomm}) can
be identified with Diagram~(\ref{spqr}) for a
suitable $A$. Let $V$ be a $6$-dimensional
complex vector-space and $F\hra\wedge^3
V\otimes\cO_{\PP(V)}$ be the sub-vector-bundle
defined by~(\ref{fibradieffe}).
\begin{prp}\label{effeomega}
Keep notation as above. Then
\begin{equation}\label{bellanapoli}
 F\cong\Omega^3_{\PP(V)}(3).
\end{equation}
\end{prp}
\begin{proof}
Let $Q:=\Theta_{\PP(V)}(-1)$. Thus we have the
Euler sequence
\begin{equation}\label{succeul}
  0\to\cO_{\PP(V)}(-1)\lra V\otimes\cO_{\PP(V)}
  \lra Q\to 0
\end{equation}
and by definition  $F\cong (\wedge^2 Q)(-1)$.
The perfect pairing $\wedge^2 Q\times\wedge^3
Q\to\wedge^5 Q\cong\cO_{\PP(V)}(1)$ gives an
isomorphism $\wedge^2 Q\cong (\wedge^3
Q^{\vee})(1)$ and hence
\begin{equation}\label{altereffe}
  F\cong \wedge^3 Q^{\vee}\cong
  \Omega^3_{\PP(V)}(3).
\end{equation}
\end{proof}
In order to identify~(\ref{grancomm}) with~(\ref{spqr}) we will need  a few
properties of the vector-bundle $F\cong \Omega^3_{\PP(V)}(3)$.
\begin{prp}\label{venti}
Keep notation as above. The dual of Exact
Sequence~(\ref{effeinwedge}) defines an
isomorphism
\begin{equation}\label{sezioni}
  \wedge^3 V^{\vee}\overset{\sim}{\lra}
  H^0(F^{\vee}).
\end{equation}
\end{prp}
\begin{proof}
Exact Sequence~(\ref{succeul}) gives an exact
sequence
\begin{equation}\label{}
  0\to
  \wedge^2 V\otimes\cO_{\PP(V)}(-1)
  \lra\wedge^3 V\otimes\cO_{\PP(V)}
  \lra\wedge^3 Q\to 0.
\end{equation}
This induces an isomorphism
\begin{equation}\label{altramappa}
\wedge^3 V\overset{\sim}{\lra}H^0(\wedge^3 Q).
\end{equation}
The symplectic form on $\wedge^3 V$ gives an identification $\wedge^3
V^{\vee}\cong \wedge^3 V$ and we have $F^{\vee}\cong \wedge^3 Q$
by~(\ref{altereffe}). With these identifications the map of~(\ref{sezioni}) is
identified with the map of~(\ref{altramappa}). Thus~(\ref{sezioni}) is an
isomorphism.
\end{proof}
\begin{prp}\label{tuttesez}
Keep notation as above. Assume that $\cW$ is a
symplectic vector-bundle and that  $\mu\colon
F\to\cW$ is an injection of vector-bundles such
that $\mu(F)$ is a Lagrangian
sub-vector-bundle. Then $\mu^{\vee}\colon
\cW^{\vee}\to F^{\vee}$ induces an isomorphism
\begin{equation}\label{}
  H^0(\cW^{\vee})\overset{\sim}{\lra}
  H^0(F^{\vee}).
\end{equation}
\end{prp}
\begin{proof}
Since $\mu(F)$ is Lagrangian the symplectic
form on $\cW$ induces an isomorphism
$\cW/\mu(F)\cong F^{\vee}$. Thus we have an
exact sequence
\begin{equation}\label{}
  0\to F\mapor{\mu}\cW\mapor{}F^{\vee}\to 0.
\end{equation}
and its dual
\begin{equation}\label{autoduale}
  0\to F\mapor{}\cW^{\vee}\mapor{\mu^{\vee}}
  F^{\vee}\to 0.
\end{equation}
The above exact sequence with $\cW=\wedge^3
V\otimes\cO_{\PP(V)}$ gives
\begin{equation}\label{svanisce}
  0=h^0(F)=h^1(F)
\end{equation}
because of Proposition~(\ref{venti}). Now
consider~(\ref{autoduale}) in general: we get
that $H^0(\cW)\mapor{} H^0(F^{\vee})$ is an
isomorphism because of
Equation~(\ref{svanisce}).
\end{proof}
By Proposition~(\ref{effeomega}) we have
$\Omega^3_{|H|^{\vee}}(3)\cong F$. By
Claim~(\ref{iniezione}) and
Proposition~(\ref{tuttesez}) we have a sequence
of  isomorphisms
\begin{equation}\label{tantiso}
 \wedge^3 V^{\vee}\mapor{\sim}
 H^0(F^{\vee})\mapor{\sim}
 H^0(\wedge^3\Theta_{|H|^{\vee}}(-3))\mapor{\sim}
 H^0(\theta)^{\vee}\oplus H^0(\theta).
\end{equation}
Let
\begin{equation}\label{}
  \rho\colon H^0(\theta)\oplus
H^0(\theta)^{\vee}\mapor{\sim} \wedge^3 V
\end{equation}
be the transpose of the composition of the maps
in~(\ref{tantiso}). Then (abusing notation)
\begin{equation}\label{mapparho}
  \rho(\Omega^3_{|H|^{\vee}}(3))=F.
\end{equation}
Thus~(\ref{grancomm}) starts looking
like~(\ref{spqr}). One missing link: we have
not proved that $\rho(H^0(\theta)^{\vee})$ is
Lagrangian for the symplectic form $\sigma$ on
$\wedge^3 V$ defined in
Section~(\ref{prologo}).
\begin{prp}\label{standard}
 Let $V$ be a $6$-dimensional complex vector-space.
 Suppose that $\tau\colon\wedge^3V\times\wedge^3
  V\to\CC$ is a symplectic form such that
  $\tau|_{F_{\ell}}=0$ for every
  $\ell\in\PP(V)$. Then $\tau=c\sigma$ for a
  certain $c\in\CC^{*}$.
\end{prp}
\begin{proof}
Let $\{e_0,\ldots,e_5\}$ be a basis of $V$. Let
$\cS\subset\cP(\{0,\ldots,5\})$ be the family
of $I\subset\{0,\ldots,5\}$ of cardinality $3$.
For $\{i,j,k\}\in\cS$ with $i<j<k$ we let
$e_I:=e_i\wedge e_j\wedge e_k$. Choose an
ordering of $\cS$; then
$\{\ldots,e_I,\ldots\}_{I\in\cS}$ is basis of
$\wedge^3 V$. Let $\alpha,\beta\in\wedge^2V$.
For $i\in\{0,\ldots,5\}$ we have
$e_i\wedge\alpha,e_i\wedge\beta\in F_{e_i}$; by
our hypothesis we get that
$\tau(e_i\wedge\alpha,e_i\wedge\beta)=0$ and
hence
\begin{equation}\label{nonvuoto}
  \text{$\tau(e_I,e_J)=0$ if $I\cap J\not=\emptyset$.}
\end{equation}
Let $j\in\{0,\ldots,5\}$; then
\begin{equation}\label{}
  0=\tau((e_i+e_j)\wedge\alpha,(e_i+e_j)\wedge\beta)=
 \tau(e_i\wedge\alpha,e_j\wedge\beta)+
 \tau(e_j\wedge\alpha,e_i\wedge\beta).
\end{equation}
This implies that
\begin{equation}\label{consegno}
  sign(I,I^c)\tau(e_I,e_{I^c})=
  sign(J,J^c)\tau(e_J,e_{J^c}),\quad I,J\in\cS
\end{equation}
where $I^c,J^c$ are the complements of $I,J$ in
$\{0,\ldots,5\}$ respectively. The proposition
is an immediate consequence of~(\ref{nonvuoto})
and~(\ref{consegno}).
\end{proof}
\begin{crl}
Keep notation and hypotheses as above. Then
\begin{equation}\label{}
  \rho^{*}\sigma=c\lambda
\end{equation}
where  $\rho$ is the isomorphism
of~(\ref{mapparho}), $\sigma$ is the symplectic
form defined in Section~(\ref{prologo}), $c$ is
a non-zero constant and $\lambda$ is  the
tautological symplectic form defined
in~(\ref{formataut}).
\end{crl}
\begin{proof}
The corollary is equivalent to the equality
\begin{equation}\label{multstan}
  (\rho^{-1})^{*}\lambda=c^{-1}\sigma,\quad c\in \CC^{*}.
\end{equation}
By Claim~(\ref{iniezione})
$\rho(\Omega^3_{|H|^{\vee}}(3))$ (yes, we abuse
notation again) is a Lagrangian
sub-vector-bundle of $\wedge^3
V\otimes\cO_{\PP(V)}$ equipped with symplectic
form $(\rho^{-1})^{*}\lambda$. By
Equality~(\ref{mapparho}) we get that for any
$\ell\in\PP(V)$ the restriction of
$(\rho^{-1})^{*}\lambda$ to $F_{\ell}$ is zero;
by Proposition~(\ref{standard}) we get
that~(\ref{multstan}) holds.
\end{proof}
Let $A:=\rho(H^0(\theta)^{\vee})$. Since
$H^0(\theta)^{\vee}$ is a Lagrangian subspace
of $H^0(\theta)\oplus H^0(\theta)^{\vee}$
equipped with the symplectic form $\lambda$ the
above corollary gives that  $A\in
\mathbb{LG}(\wedge^3 V)$. It is clear from
Diagram~(\ref{grancomm}) that we have equality
of reduced $Y=(Y_A)_{red}$ where $(Y_A)_{red}$
is the reduced $Y_A$. Since $\deg
Y=6=\deg(Y_A)$ and $Y,Y_A$ are both Cartier
divisors we get that
\begin{equation}\label{yeya}
  Y=Y_A.
\end{equation}
\begin{clm}
Keep notation and assumptions as above.
Then $A\in \mathbb{LG}(\wedge^3 V)^0$.
\end{clm}
\begin{proof}
First notice that $A\in \mathbb{LG}(\wedge^3 V)^{\times}$ simply because
$j_{*}\theta$ is a Casnati-Catanese sheaf. We notice also that
\begin{equation}\label{baba}
  D_1(A,F)\setminus D_2(A,F)=Y\setminus sing Y.
\end{equation}
Thus we may apply Proposition~(\ref{yewgen}) in order to prove the claim.
Item~(3) of the proposition is satisfied by~(\ref{baba}) and~(\ref{eqloc}) and
hence $A\in \mathbb{LG}(\wedge^3 V)^0$. Alternatively one can check that
Item~(3) of Proposition~(\ref{yewgen}) holds. The only unproved fact is that
$W_A$ is smooth. By~(\ref{diconipull})-(\ref{codimsigma}) and~(\ref{pusingy}) we
know that $W_A$ is a local complete intersection; since $sing Y=(W_A)_{red}$ and
$sing Y$ is smooth we get that if $W_A$ is not smooth then it is not reduced.
Thus by Formula~(\ref{classew}) we get that it suffices to show that $\deg(sing
Y)=40$; this follows at once form the formulae in Item~(1) of Theorem~(1.1)
of~\cite{ognum}.
\end{proof}
Comparing Diagrams~(\ref{spqr})
and~(\ref{grancomm}) we see that
$\theta=\zeta_A$; by~(\ref{defxia})
and~(\ref{tetadef}) we get that $\eta=\xi_A$.
This completes the proof of
Theorem~(\ref{chiappone}).
\section{An involution on a moduli space}
 \label{involsuk}
 \setcounter{equation}{0}
 Let $\cK^0_2$ be the set of isomorphism classes
 of couples
 $(X,H)$ where $X$ is a
 numerical $(K3)^{[2]}$ and $H$ an ample
 divisor on $X$ such that~(\ref{gradodod}) and~(a)
 of Section~(\ref{prologo}) both hold; couples
 $(X_i,H_i)$ for $i=1,2$ are isomorphic if
 there exists an isomorphism $\psi\colon
 X_1\overset{\sim}{\lra} X_2$ such that
 $\psi^{*}H_2\sim H_1$. By
 Theorem~(\ref{mainthm}) we have an
 identification
\begin{equation}\label{}
  \cK^0_2=
 \mathbb{LG}(\wedge^3 V)^0//\mathbb{PGL}(V)
\end{equation}
where $V$ is a $6$-dimensional complex
vector-space. We remark that the
 second equality
of~(\ref{gradodod}) follows from the first one.
Furthermore the first equality
of~(\ref{gradodod}) should be thought of as the
analogue of self-intersection $2$ for an ample
divisor on a $K3$ surface, see Footnote~(1). We
recall also that Condition~(a) is an open
condition. Thus $\cK^0_2$ is an open subset of
the moduli space of couples $(X,H)$ where $X$ is
a numerical $(K3)^{[2]}$ and $H$ is an ample
divisor on $X$ of square $2$ for Beauville's
quadratic form - actually the larger moduli space
of couples $(X,H)$ with $H$ big and nef is a
better setting for what follows. We represent a
point of $\cK^0_2$ by $[X_A]$ where $A\in
\mathbb{LG}(\wedge^3 V)^0$. Let $\delta_V\colon
\mathbb{LG}(\wedge^3
V)\overset{\sim}{\lra}\mathbb{LG}(\wedge^3
V^{\vee})$ be the isomorphism
of~(\ref{annichilo}). The subset of
$\mathbb{LG}(\wedge^3 V)^0$ appearing
in~(\ref{pigna}) is open dense and
$\mathbb{PGL}(V)$-invariant hence
\begin{equation}\label{}
  \cU:=\left(\mathbb{LG}(\wedge^3 V)^0\cap
  \delta_V^{-1}\mathbb{LG}(\wedge^3 V^{\vee})^0\right)
 //\mathbb{PGL}(V).
\end{equation}
 is an open and dense subset of $\cK^0_2$. If
$[X_A]\in\cU$ then we have the nautural double
cover $X_{A^{\bot}}\to Y_{A^{\bot}}$ and
furthermore $[X_{A^{\bot}}]\in\cK^0_2$ by
Theorem~(\ref{mainthm}). Hence duality defines an
involution
\begin{equation}\label{}
  \begin{matrix}
 \cU & \mapor{\delta} & \cU \\
 [X_A] & \mapsto & [X_{A^{\bot}}]
  \end{matrix}
\end{equation}
Let's prove that $\delta$ is not the identity. We
consider the example of Mukai that we presented
in the proof that Theorem~(\ref{chiappone})
implies Theorem~(\ref{mainthm}). Thus $S$ is a
$K3$ given by~(\ref{kappatre})
satisfying~(\ref{dieci}), and $Y,X$ are given
by~(\ref{eccoy}) and~(\ref{eccox}) respectively.
Then $Y\subset |I_S(2)|$ corresponds to a certain
$M\in \mathbb{LG}(\wedge^3 H^0(I_S(2)))^0$
i.e.~$Y=Y_M$. According to Subsubsection~(5.4.2)
of~\cite{oginv} the dual $Y_M^{\vee}\subset
|I_S(2)|^{\vee}$ is described as follows.
By~(\ref{dieci}) $S$ is cut out by quadrics and
contains no lines hence we have a well-defined
regular map
\begin{equation}\label{fdoppiav}
 \begin{matrix}
S^{[2]} & \mapor{g} & |I_S(2)|^{\vee}\\
[Z] & \mapsto & \{Q|\ Q\supset\langle
Z\rangle\}
\end{matrix}
\end{equation}
where $\langle Z\rangle\subset\PP^6$ is the
unique line containing $Z$. Then
$Y_M^{\vee}=Im(g)$. One has $|I_F(2)|\in Im(g)$ -
in fact $|I_F(2)|$ is the image by $g$ of
\begin{equation}\label{}
  B_S:=\{[Z]\in S^{[2]}|\ \langle Z\rangle
  \subset F\}.
\end{equation}
Iskovskih (Cor.~(6.6) of~\cite{isk}) proved that
$B_S\cong\PP^2$. We proved that $g$ has degree
$2$ onto its image and that
\begin{equation}\label{gradoduale}
  \deg Y_M^{\vee}=6.
\end{equation}
For $[Z]\in(S^{[2]}\setminus B_S)$ there exists a
unique conic $C\subset F$ containing $Z$. Then
$C\cap S$ is a scheme of length $4$ containing
$Z$ and there is a well-defined residual scheme
$Z'$ of $Z$, see Subsection~(4.3)
of~\cite{oginv}. Let $\phi([Z]):=[Z']$; then
\begin{equation}\label{fibradig}
  g^{-1}(g([Z]))=\{[Z],\phi([Z])\}.
\end{equation}
For $[Z]$ generic $Z'$ can be characterized as
the unique $Z'\subset S$ of length $2$ such that
\begin{equation}\label{}
 Z'\cap Z=\emptyset,\quad
  \langle Z'\rangle \cap \langle Z\rangle
  \not=\emptyset
\end{equation}
\begin{prp}\label{singdiy}
Keep notation and assumptions as above. Then
$|I_F(2)|\in Y_M^{\vee}$ is a point of
multiplicity $3$.
\end{prp}
\begin{proof}
 Let $\Lambda\subset |I_F(2)|$ be a generic
 linear subspace with $\dim\Lambda=3$. Thus
\begin{equation}\label{}
  \Lambda^{\bot}:=\{\Omega\in |I_S(2)|^{\vee}\,|\
  \Lambda\subset\Omega\}
\end{equation}
is a generic line in $|I_S(2)|^{\vee}$ containing
$|I_F(2)|$. By~(\ref{gradoduale}) we must prove
that
\begin{equation}\label{sonotre}
  \sharp\left(\Lambda^{\bot}\cap Y_M^{\vee}
  \setminus\{|I_F(2)|\}\right)=3.
\end{equation}
Let $[Z]\in S^{[2]}$; then
$g([Z])\in\left(\Lambda^{\bot}\cap Y_M^{\vee}
  \setminus\{|I_F(2)|\}\right)$ if and only if
\begin{equation}\label{dentrofuori}
  \langle Z\rangle\subset\bigcap
  \limits_{t\in\Lambda}Q_t
  \quad
 \langle Z\rangle\not\subset F.
\end{equation}
(Here $Q_t\subset\PP^6$ is the quadric
corresponding to $t$.) In the proof of Item~(1)
of Lemma~(4.20) of~\cite{oginv} we showed that
\begin{equation}\label{}
  \bigcap
  \limits_{t\in\Lambda}Q_t=F\cup\langle
  C\rangle
\end{equation}
where $C\subset F$ is a certain conic with span
$\langle C\rangle$ such that $F\cap\langle
C\rangle=C$. Let $\ov{Q}\subset\PP^6$ be a
quadric such that $S=F\cap\ov{Q}$,
see~(\ref{kappatre}). By~(\ref{dieci}) the
surface $S$ does not contain conics and hence
$\ov{Q}\cap C$ is a finite set of length $4$.
Since $\Lambda$ is chosen generically
$\ov{Q}\cap C$ consists of $4$ distinct points
$p_1,\ldots,p_4$ and hence
\begin{equation}\label{}
 \langle C\rangle\cap S=\{p_1,\ldots,p_4\}.
\end{equation}
Thus $[Z]\in S^{[2]}$
satisfies~(\ref{dentrofuori}) if and only if
$Z$ is a subset of $\{p_1,\ldots,p_4\}$; since
there are $6$ such $Z$'s and since $g$ has
degree $2$ onto its image
(see~(\ref{fibradig})) we get
that~(\ref{sonotre}) holds.
\end{proof}
The above proposition shows that the involution
$\delta$ is not the identity. In fact $M\in
\mathbb{LG}(\wedge^3 H^0(I_S(2)))^0$ and hence
every singular point of $Y_M$ has multiplicity
$2$. Since $Y_M^{\vee}$ has a point of
multiplicity $3$ we get that $Y_M^{\vee}$ is not
projectively isomorphic to $Y_M$. The proposition
also shows that if we want the codomain of
$\delta$ to be $\cK^0_2$ then $\delta$ is not
defined at $Y_M$. Let $\cK_2$ be the set of
isomorphism classes
 of couples $(X,H)$ where $X$ is a
 numerical $(K3)^{[2]}$, $H$ is a big and nef
 divisor on $X$ such that~(\ref{gradodod})
 holds and furthermore $(X,H)$ is a \lq\lq
 limit\rq\rq of $(X',H')$ parametrized by
 $\cK^0_2$; is it true that $\delta$ extends to
 a regular involution defined on all of
 $\cK_2$?
 \section{Particular examples of EPW-sextics}
 \label{esespl}
 \setcounter{equation}{0}
 The question we briefly address is the
 following. Given $(X,H)$
 satisfying~(\ref{gradodod}) and Item~(a) of
 Section~(\ref{prologo}) how do we describe an $A\in
 \mathbb{LG}(\wedge^3 H^0(\cO_X(H))^{\vee})$
 such that $X_A\cong X$? More generally we may
 ask the same question for an $(X',H')$ which is
 a limit of $(X,H)$ as above. A similar question
 for the Fano variety of lines on a cubic
 $4$-fold is studied by Hassett
 in~\cite{hassett}.
\subsection{Another description of $Y_A$}
 Let $V$ be a $6$-dimensional complex vector-space.
 Let $\mathbb{LG}(\wedge^3 V)^{\dag}
 \subset\mathbb{LG}(\wedge^3 V)$ be the open subset of
 $A$ such that there exists $W\subset V$ of
codimension $1$ such that
\begin{equation}\label{}
  \wedge^3 W\cap A=\emptyset.
\end{equation}
For $A\in\mathbb{LG}(\wedge^3 V)^{\dag}$ we will
describe $Y_A$
 as a set of degenerate quadrics in $\PP(\wedge^2 W)$.
 Using this description
 we will propose an answer to the question asked
 above for those $Y$ that were
 described by Mukai in Ex.(5.17) of~~\cite{mukai}
  and that we
 have used in
 Sections~(\ref{dimprinc})-(\ref{involsuk}).
Let $\ell_0\in\PP(V)$ be such that
\begin{equation}\label{}
  V=\ell_0\oplus W.
\end{equation}
Thus we have a decomposition
\begin{equation}\label{partictriv}
  \wedge^3 V=\left(\ell_0\otimes\wedge^2 W\right)
  \oplus\wedge^3 W.
\end{equation}
Both addends are Lagrangian subspaces of
$\wedge^3 V$ and hence the symplectic form
$\sigma$ induces an isomorphism
\begin{equation}\label{sonoduali}
  \wedge^3 W\overset{\sim}{\lra}\left(\ell_0\otimes\wedge^2
  W\right)^{\vee}.
\end{equation}
We choose a non-zero $v_0\in\ell_0$.
Multiplication by $v_0$ defines an isomorphism
$\wedge^2 W \to \ell_0\otimes\wedge^2 W$ and
hence~(\ref{sonoduali}) becomes an isomorphism
\begin{equation}\label{soduali}
  \wedge^3 W\overset{\sim}{\lra}\wedge^2
  W^{\vee}.
\end{equation}
Let $vol_W\colon\wedge^5
W\overset{\sim}{\lra}\CC$ be the trivialization
defined by setting
$vol_W(\tau)=vol(v_0\wedge\tau)$; then
Isomorphism~(\ref{soduali}) is given by the
perfect pairing
\begin{equation}\label{}
   \begin{matrix}
   \wedge^2 W\times\wedge^3 W & \lra & \CC  \\
 (\alpha,\beta) & \mapsto & vol_W(\alpha\wedge\beta)
 \end{matrix}
\end{equation}
Let
\begin{equation}\label{}
  U:=\PP(V)\setminus\PP(W).
\end{equation}
Tensorizing~(\ref{partictriv}) by $\cO_U$ we get
a symplectic trivialization of $\wedge^3
V\otimes\cO_U$ with $\cL:=\wedge^2 W\otimes\cO_U$
and $\cL^{\vee}:=\wedge^3 W\otimes\cO_U$. Notice
that $F_{\ell}$ is transversal to $\wedge^3 W$
for every $\ell\in U$. Thus
$D_1(A,F,U,\cL,\cL^{\vee})$ is well-defined. Let
$q_A\in Sym^2(\wedge^2 W^{\vee})$ be as
in~(\ref{mappegrafi}). For $w\in W$ let $q_w\in
Sym^2(\wedge^2 W^{\vee})$ be the quadratic form
defined by
\begin{equation}\label{}
  q_w(\alpha):=w\wedge\alpha\wedge\alpha.
\end{equation}
Let $\ell=[v_0+w]$; as is easily checked $q_w$ is
the quadratic form defined by the symmetric map
$\phi_{F_{\ell}}$. The map
\begin{equation}\label{}
  \begin{matrix}
 W & \lra & U \\
 w & \mapsto & [v_0+w]
  \end{matrix}
\end{equation}
is an isomorphism; from now on we identify $U$
with $W$ via the above map. Thus
\begin{equation}\label{}
  D_1(A,F,U,\cL,\cL^{\vee})=Y_A\setminus\PP(W)=
  \{w\in W|\ \det(q_A-q_w)=0\}.
\end{equation}
Now  notice that the $q_w$'s are the Pl\"ucker
quadratic forms whose zero-locus is ${\bf
Gr}(2,W)\subset\PP(\wedge^2 W)$. Let $Q_A\subset
\PP(\wedge^2 W)$ be $Q_A=V(q_A)$ and let
$Z_A\subset \PP(\wedge^2 W)$ be
\begin{equation}\label{}
  Z_A:=Q_A\cap{\bf Gr}(2,W).
\end{equation}
Thus $|I_{Z_A}(2)|$ is the span of $Q_A$ and
$|I_{{\bf Gr}(2,W)}(2)|$, in particular
$|I_{Z_A}(2)|\cong\PP^5$. Let $\Sigma_A$ be the
degree-$10$ divisor on $|I_{Z_A}(2)|$
 parametrizing singular quadrics.
 Each quadric $V(q_w)\in |I_{{\bf Gr}(2,W)}(2)|$
 is singular with $\dim(sing V(q_w))=3$ and
 hence we have
\begin{equation}\label{}
  \Sigma_A=4|I_{{\bf Gr}(2,W)}(2)|+\Sigma'_A.
\end{equation}
Thus $\Sigma'_A$ is a degree-$6$ divisor. The
above discussion proves the following result.
\begin{prp}
Keep notation and assumptions as above. Then
\begin{equation}\label{}
  Y_A\cong\Sigma'_A.
\end{equation}
\end{prp}
Now we can propose an explicit description of
those $M\in\mathbb{LG}(\wedge^3 V)^0$ such that
$Y_M$ is one of Mukai's examples. Actually we
propose a description of those
$A\in\mathbb{LG}(\wedge^3 V)$ for which
$Y_A^{\vee}\cong Y_M$; if the description is
correct then one simply sets $M=A^{\bot}$. By
Proposition~(\ref{singdiy}) we know that
$Y_M^{\vee}$ has a point of multiplicity $3$
hence we choose $A$ such that
\begin{equation}\label{singqa}
  \dim(sing Q_A)=2
\end{equation}
because $\Sigma'_A$ will have a point of
multiplicity at least $3$ at $Q_A$ - of
multiplicity equal to $3$ if $A$ is generic. The
question is: is it true that if~(\ref{singqa})
holds then
\begin{equation}\label{ecchime}
  Y_A\cong g(S^{[2]})
\end{equation}
where $S \subset\PP^6$ is a certain $K3$ surface
of genus $6$ and $g$ is as in
Section~(\ref{involsuk})? Our observation is that
one can associate to such an $A$ a $K3$ surface
$S\subset\PP^6$. In fact consider the duals
$Q_A^{\vee},{\bf Gr}(2,W)^{\vee}\subset
\PP(\wedge^2 W)^{\vee}$. By~(\ref{singqa})
$Q_A^{\vee}$ is a smooth quadric hypersurface in
$sing(Q_A)^{\bot}\cong\PP^6$ and ${\bf
Gr}(2,W)^{\vee}={\bf Gr}(2,W^{\vee})$. Thus
\begin{equation}\label{}
  S:=Q_A^{\vee}\cap{\bf
Gr}(2,W)^{\vee}
\end{equation}
 is indeed a $K3$ surface of genus $6$; our guess
 is that~(\ref{ecchime}) holds with the above $S$.
\subsection{Non-reduced $Y_A$'s}
We will consider examples of  $(X_0,H_0)$ a limit
of $(X,H)$ with $X$ a numerical $(K3)^{[2]}$ and
$H$ an ample divisor on $X$ satisfying
both~(\ref{gradodod}) and Item~(a) of
Section~(\ref{prologo}) and such that $|H_0|$ is
base-point free with
\begin{equation}\label{eccemappa}
  X_0\lra|H_0|^{\vee}
\end{equation}
of degree higher than $2$, namely $4$ or $6$. In
the examples  there is an anti-symplectic
involution $\phi_0\colon X_0\to X_0$ with
quotient map $f_0\colon X_0\to Y_0$. There is a
map $j_0\colon Y_0\to |H_0|^{\vee}$ such that
Map~(\ref{eccemappa}) is the composition
\begin{equation}\label{}
  X_0\mapor{f_0}Y_0\mapor{j_0}|H_0|^{\vee}.
\end{equation}
Of course $j_0$ is not an emebedding, it is
finite onto its image. We consider the
decomposition
$f_{0,*}\cO_{X_0}=\cO_{Y_0}\oplus\eta$ where
$\eta$ is the $(-1)$-eigenspace for the action of
$\phi_0^{*}$. Let $\theta:=\eta\otimes
j_0^{*}\cO_{|H_0|^{\vee}}(3)$. Assuming that
$\theta$ is globally generated we may consider
the exact sequence
\begin{equation}\label{}
  0\to G\mapor{}
  H^0(\theta)\otimes\cO_{|H_0|^{\vee}}
  \mapor{\epsilon}j_{0,*}\theta\to 0
\end{equation}
where $\epsilon$ is the evaluation map. Then
$G\cong \Omega^3(3)_{|H_0|^{\vee}}$. One can
construct a commutative diagram as
in~(\ref{grancomm}) and proceed as in
Section~(\ref{dimprinc}) to get a decomposition
into Lagrangian subspaces
\begin{equation}\label{prego}
 \wedge^3 H^0(\cO_{X_0}(H_0))^{\vee}=H^0(\theta)
 \oplus H^0(\theta)^{\vee}.
\end{equation}
Then  $H^0(\theta)^{\vee} \in\mathbb{LG}(\wedge^3
H^0(\cO_{X_0}(H_0))^{\vee})$ is a point
corresponding to $X_0$. The question of course is
to describe Decomposition~(\ref{prego}). Our
first example is $X_0=S^{(2)}$ where $\pi\colon
S\to\PP^2$ is a double cover branched over a
smooth sextic. Thus $S$ is a $K3$ surface and
$H_S:=\pi^{*}\cO_{\PP^2}(1)$ is an ample divisor
on $S$ with $H_S\cdot H_S=2$. Let $L\subset\PP^2$
be a line and let
\begin{equation}\label{}
  D_L:=\{p_1+p_2\in S^{(2)}|\
  \{p_1,p_2\}\cap L\not=\emptyset\}.
\end{equation}
Then $D_L$ is an ample Cartier divisor on $X_0$,
call it $H_0$ (thus the pull-back of $H_0$ by the
desingularization $S^{[2]}\to S^{(2)}$ is big and
nef). There exist  smoothings $X$ of $X_0$ for
which $H_0$ deforms to an ample divisor $H$ on
$X$. Then $X$ is a deformation of $(K3)^{[2]}$
and $H$ satisfies~(\ref{gradodod}) and Item~(a)
of Section~(\ref{prologo}). The involution
$\phi_0\colon X_0\to X_0$ is given by
\begin{equation}\label{involtino}
  \phi_0(p_1+p_2)=\iota(p_1)+\iota(p_2)
\end{equation}
where $\iota\colon S\to S$ is the covering
involution of $\pi$. Let $q\colon S^2\to S^{(2)}$
be the quotient map; then
$q^{*}\cO_{X_0}(H_0)\cong\cO_{S}(H_S)\boxtimes
\cO_{S}(H_S)$ and hence
\begin{equation}\label{tiro}
 q^{*} H^0(\cO_{X_0}(k H_0))=
 Sym^2 H^0(\cO_{S}(k H_S)).
\end{equation}
Let $U:=H^0(\cO_S(H_S))=H^0(\cO_{\PP^2}(1))$.
From~(\ref{tiro}) we get that
$H^0(\cO_{X_0}(H_0))\cong Sym^2 U$. Furthermore
the image of $j_0\colon Y_0\to \PP(Sym^2
U^{\vee})$ is the chordal variety of the Veronese
surface, i.e.~the discriminant cubic, and $j_0$
is of degree $2$ onto its image. Let $\sigma\in
H^0(\cO_S(3H_S))$ be a generator of the
$(-1)$-eigenspace for the action of $\iota^{*}$;
thus the divisor of $\sigma$ is the ramification
divisor of $\pi$. One easily checks that we have
an isomorphism
\begin{equation}\label{vincenzone}
  \begin{matrix}
 Sym^3 U & \mapor{\sim} & q^{*}H^0(\theta) \\
 \alpha & \mapsto &
 \sigma\otimes\pi^{*}\alpha+\pi^{*}\alpha\otimes\sigma
  \end{matrix}
\end{equation}
(Here $\pi^{*}\colon Sym^3
U=H^0(\cO_{\PP^2}(3))\hra H^0(\cO_S(3H_S))$.) On
the other hand we have the decomposition of
$GL(U)$-modules
\begin{equation}\label{}
  \wedge^3 H^0(\cO_{X_0}(H_0))^{\vee}=
  \wedge^3(Sym^2 U^{\vee})=Sym^3 U\oplus
 Sym^3 U^{\vee}.
\end{equation}
One can check that given
Isomorphism~(\ref{vincenzone}) this is
Decomposition~(\ref{prego}) in the present case.
 The last example is similar but we must state that we
  have not
 checked the details. Let $S\subset\PP^3$ be a
 smooth quartic containing no lines and let
 $X_0:=S^{[2]}$. Let $U:=H^0(\cO_S(1))$.
 We have a map
\begin{equation}\label{}
   \begin{matrix}
 X_0=S^{[2]} & \mapor{g} & {\bf Gr}(2,U^{\vee})\\
 [Z] & \mapsto & \langle Z\rangle
   \end{matrix}
\end{equation}
Let $p\colon {\bf Gr}(2,U^{\vee})\hra\PP(\wedge^2
U^{\vee})=\PP^5$ be the Pl\"ucker embedding and
let $H_0:=(pg)^{*}\cO_{\PP^5}(1)$ be the
Pl\"ucker class. Then $(X_0,H_0)$ is the limit of
$(X,H)$ where $H$ satisfies~(\ref{gradodod})
and~(a) of Section~(\ref{prologo}). The
involution $\phi_0\colon X_0\to X_0$ associates
to $[Z]$ the residual of $Z$ in $\langle Z\rangle
\cap S$. Then $j_0\colon Y_0\to {\bf
Gr}(2,U^{\vee})$ is finite of degree $3$. Since
$H^0(\cO_{X_0}(H_0))=\wedge^2 U$ it is natural to
guess that in the present case
Decomposition~(\ref{prego}) is the decomposition
of $GL(U)$-modules
\begin{equation}\label{}
  \wedge^3(\wedge^2 U^{\vee})=
  Sym^2 U\oplus Sym^2 U^{\vee}.
\end{equation}
\vskip 1cm
 \scriptsize{
Kieran G. O'Grady\\
Universit\`a di Roma ``La Sapienza",\\
Dipartimento di Matematica ``Guido Castelnuovo",\\
Piazzale Aldo Moro n.~5, 00185 Rome, Italy,\\
e-mail: {\tt ogrady@mat.uniroma1.it}. }

\end{document}